\documentclass{article} 
\usepackage{amssymb,amsmath,amsthm,graphicx,geometry,comment,xcolor,mathtools,wrapfig}
\usepackage[margin=2cm,font=scriptsize]{caption}

\usepackage[pagebackref,breaklinks=true]{hyperref}\hypersetup{colorlinks,
  linkcolor={green!50!black},
  citecolor={green!50!black},
  urlcolor=blue
}

\newcommand{\pp}{/\hspace{-3pt}/}
\newcommand{\p}{\partial}

\newcommand\DG{\reflectbox{\rotatebox[origin=c]{180}{$\mathbb L$}}}

\title{Polynomial Invariant of Molecular Circuit Topology}

\author{Alireza Mashaghi\footnote{Leiden Academic Centre for Drug Research, Faculty of Mathematics and Natural Sciences, Leiden University, 2311 Leiden, The Netherlands} and Roland van der Veen\footnote{Bernoulli Institute, Department of Mathematics, Faculty of Science and Engineering, University of Groningen, 9712 Groningen, The Netherlands; r.i.van.der.veen@rug.nl}} 
\begin{document}

\abstract{The topological framework of circuit topology has recently been introduced to complement knot theory and to help in understanding the physics of molecular folding. Naturally evolved linear molecular chains, such as proteins and nucleic acids, often fold into 3D conformations with critical chain entanglements and local or global structural symmetries stabilised by formation contacts between different parts of the chain. Circuit topology captures the arrangements of intra-chain contacts within a given folded linear chain and allows for the classification and comparison of chains. Contacts keep chain segments in physical proximity and can be either mechanically hard attachments or soft entanglements that constrain a physical chain. Contrary to knot theory, which offers many established knot invariants, circuit invariants are just being developed. Here, we present polynomial invariants that are both efficient and sufficiently powerful to deal with any combination of soft and hard contacts. A computer implementation and table of chains with up to three contacts is also~provided.}

\keyword{circuit topology; knot theory; invariant; folding; polymer} 

\section{Introduction}

Linear polymers are an important subset of macromolecules with critical roles in living organisms and are used in engineering applications~\cite{origami20}. A~linear polymer is a molecular chain made of units, so called monomers. By~changing the chemical properties of these monomers, one can make polymers with different physicochemical properties. This fundamental concept in chemistry has led to emergence of many synthetic polymers with various applications in medicine and industry. Living organisms however generate a wide range of linear polymers using a limited set of monomer chemistries. This is because biomolecular chains, such as proteins and nucleic acids, typically fold into 3D conformations which give the molecules new properties. Proteins can form various folded structures at various scales and with different symmetries by forming intra-chain contacts; proteins may also form knots and slipknots~\cite{JS15,FHW19,origami20}. Inspired by biology, chemists have only recently started synthesizing folded molecular chains~\cite{tezuka19}. Molecular engineering typically uses a bottom-up approach which involves synthesizing basic fold units and then connecting them to generate complexity. Generating complex folded linear chains requires advancements in our synthetic methodology as well as suitable conceptual mathematical framework for characterization and comparison of topological complexity. The~latter inspired the development of molecular circuit topology, a~framework that categorizes the arrangement of contacts in a folded linear chain~\cite{Ma21, MWT14}. 

Circuit topology is inspired by physics of polymers and recognises that contacts restrain dynamics of a chain and keep segments in close proximity~\cite{Heidari20}. Topological arrangement of the contacts is closely related to kinetics of their formation~\cite{Heidari20, Mugler14}. Furthermore, chain entanglement may also restrain a physical chain and effectively stabilise certain folds. Let us consider a folded linear chain. To~describe a folded chain in space there are roughly two aspects to address: First, intra-chain contacts or bonds  that turn the chain into a special type of graph. Second, the~chain sits in three dimensional space in a certain way, allowing it to form knots and tangles. Both the bonding and the tangling can constrain the chain. The~first type, we call hard contacts (or H-contacts), while the second type we will describe in terms of soft contacts (or S-contacts). Circuit topology uses a uniform language to categorise the arrangement of hard and soft~contacts.

\begin{figure}[H]
\includegraphics[width=\linewidth]{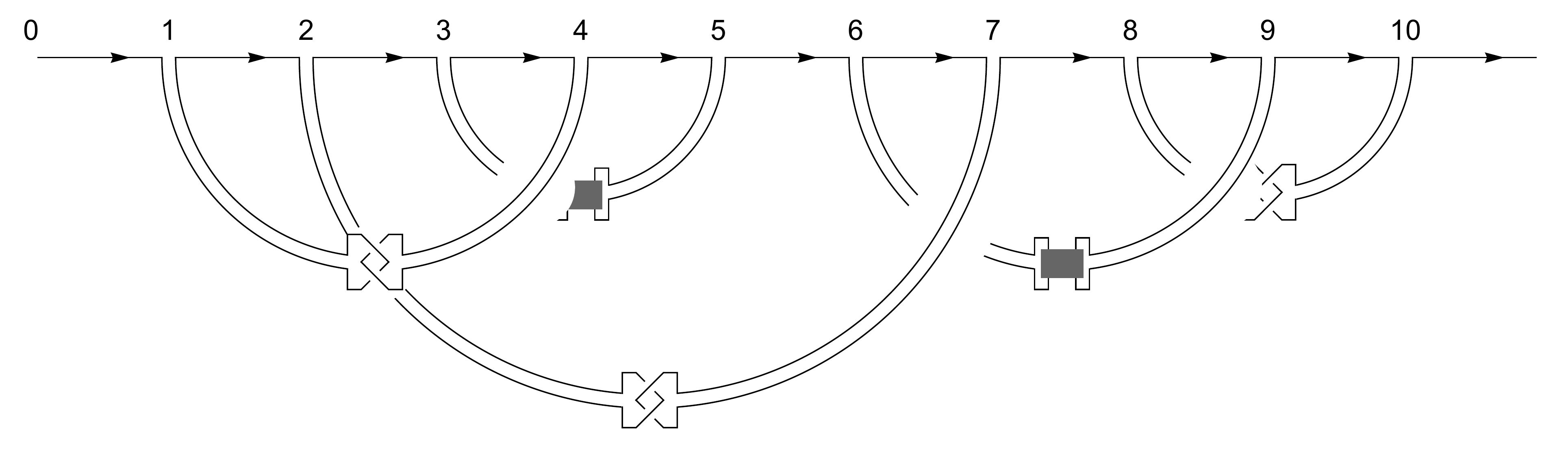}
\caption{A knotted folded chain with three soft contacts and two hard~contacts.}
\label{fig.gCTExample1}
\end{figure}

In this article, we propose a precise mathematical model for folded chains called \mbox{$H$-tangle} diagrams.  Tangles are a commonly used bottom-up approach to knot theory~\cite{Co70} and in this work we extend this approach to include hard contacts. An~example of a \mbox{$H$-tangle} appears in Figure~\ref{fig.gCTExample1}. $H$-tangles take into account both hard contacts (shown in black) and soft contacts where the chain constrains itself by a clasp. 
We describe an algorithm to assign to each $H$-tangle diagram $D$ a certain polynomial called $\DG_1(D)$.
For example $\DG_1$ applied to the example in Figure~\ref{fig.gCTExample1} is
\begin{align*}
&h_1 t^6-2 h_1 t^5-h_1 t^4+6 h_1 t^3-h_1 h_2 t^3+h_2 t^3-6 h_1 t^2-2 h_2 t^2-2 h_1t^{-2}+10 h_1 h_2t^{-2}+\\
&2 h_2t^{-2}+h_1t^{-3}+13 h_2t^{-3}-4 h_1 h_2t^{-4}-11 h_2t^{-4}+   h_1 h_2t^{-5}-3 h_2t^{-5}+5 h_2t^{-6}
   -h_2t^{-7}+\\&-3 h_1 t+5 h_1 h_2 t+h_2 t
   -2h_1t^{-1}-5 h_1 h_2t^{-1}-11 h_2t^{-1}+7 h_1-5 h_1 h_2+5 h_2+\\&-t^6+2 t^5-2 t^4+7 t^2-8t^{-2}+3t^{-4}-t^{-5}-10 t+9t^{-1}+2.
\end{align*}

The polynomial $\DG_1$ is a topological invariant of the folded chain in the sense that if $H$-tangle diagrams $D$ and $D'$ represent the same folded chain, then $\DG_1(D) = \DG_1(D')$. In~\cite{CEM21,GM20} other invariants of folded chains were proposed. The~difference is that $\DG_1$ is easier to compute in practice while still being sufficiently powerful to distinguish many different configurations. A~computer implementation of our algorithm appears in the~Appendix~\ref{sec.Appendix}.

Even though folded chains can take many forms, we can always bring them in a form such as that in Figure~\ref{fig.gCTExample1}
called a generalized circuit topology diagram~\cite{GM20}.
The chain runs horizontally from left to right, except for making a detour to make a soft or hard contact with another piece of the chain in a controlled way.
In the case of knots a similar form appears in~\cite{MP19}.
In the Appendix~\ref{sec.Appendix} we list all circuit diagrams with at most three contacts together with their corresponding value of $\DG_1$.~There do exist several circuit topology diagrams for the same folded chain
but in combination with our invariant this provides a first step towards classification of folded chains of low~complexity.

The plan of the paper is to first describe some basic notions in circuit topology using hard contacts only. 
Historically, circuit topology was first introduced to classify such circuits. Next we review the knot theory techniques that are fundamental to handling soft contacts, this time ignoring hard contacts.
We set up the Alexander polynomial invariant for knots and tangles in the form of $\Gamma$ calculus. To~include hard contacts in our tangle language we develop a theory of $H$-tangles that is similar to that of singular knots. The~invariant $\DG_1$ is then introduced as a suitable extension of $\Gamma$ calculus to $H$-tangles. Finally we present a special class of $H$-tangle diagrams especially suited for discussion of circuit topology and tabulate $\DG_1$ on~those.

\section{Circuit~Topology}

\begin{figure}[H]
\includegraphics[width=\linewidth]{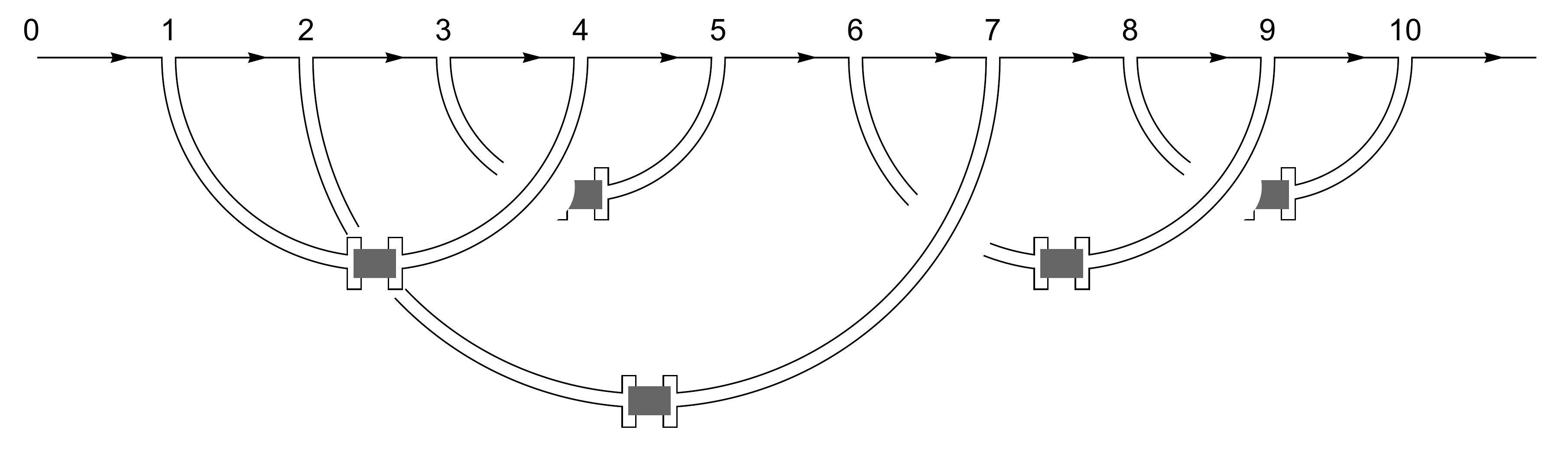}
\caption{The chain $\mathcal{F}$ with five hard~contacts.}
\label{fig.CTExample}
\end{figure}

Perhaps the most fundamental aspect of a folded chain is the way it bonds to itself. In~this section, we focus on this aspect only, ignoring for now how it sits in space. In~other words we study the folded chain as an abstract directed graph, see for example the chain $\mathcal{F}$ in Figure~\ref{fig.CTExample}.

To describe such chains in general we enumerate the spots where the chain bonds to itself in order of appearance as one walks along the chain. Assuming each spot bonds with precisely one other spot there must be an even number of such spots that we can number $1,2,\dots, 2n$. In~the example $\mathcal{F}$, we chose $n=5$ as there are five bonds. Each bond is described by the pair of integers $\{i,j\}$ indicating the two locations involved. We often refer to such bonds as hard contacts when we want to emphasize the distinction with the knotting behaviour described in later~sections.

Any folded chain is then described precisely by a partition of the integers $1,\dots 2n$ into $n$ disjoint pairs. For~example we describe $\mathcal{F}$ by 
\[\mathcal{F}=\{\{1,4\},\{2,7\},\{3,5\},\{6,9\},\{8,10\}\}\].

In listing such examples our convention is to order the list of pairs lexicographically and refer to any bond by the spot where it first appears as one walks along the chain. In~the mathematical literature such structures appear under the name \emph{chord diagrams} \cite{CDM12}.

At this point, we should emphasize that, contrary to the picture shown, the~description $\mathcal{F}$ is completely abstract. In~the picture the bond $\{1,4\}$ appears to pass in front of bond $\{2,7\}$ but this is not part of our mathematical description of the graph. So far the space around the graph has not been taken into account so there is no notion of behind or in front. Modelling and describing how abstract folded chains interact with ambient space is the subject of the next sections as it is considerably more~difficult.

\begin{figure}[H]
\includegraphics[width=5.8cm]{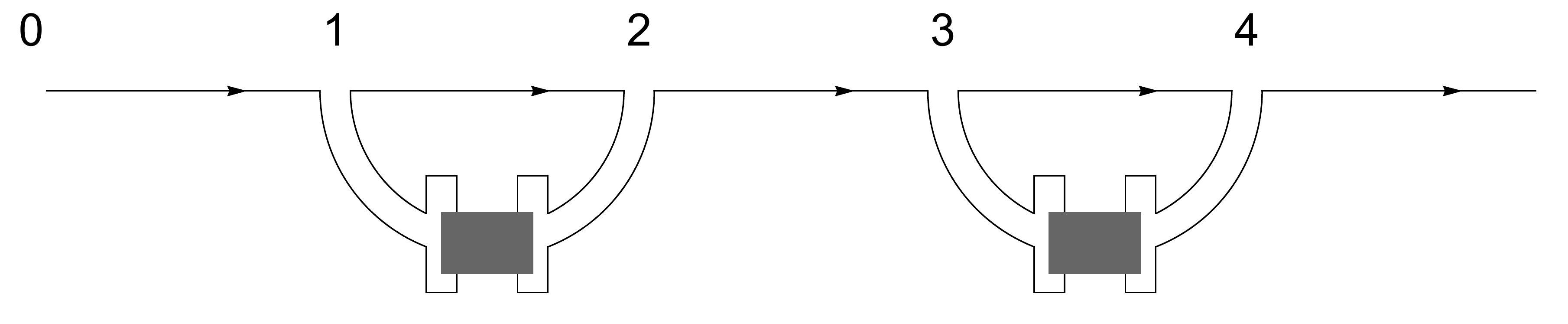}
\includegraphics[width=5.8cm]{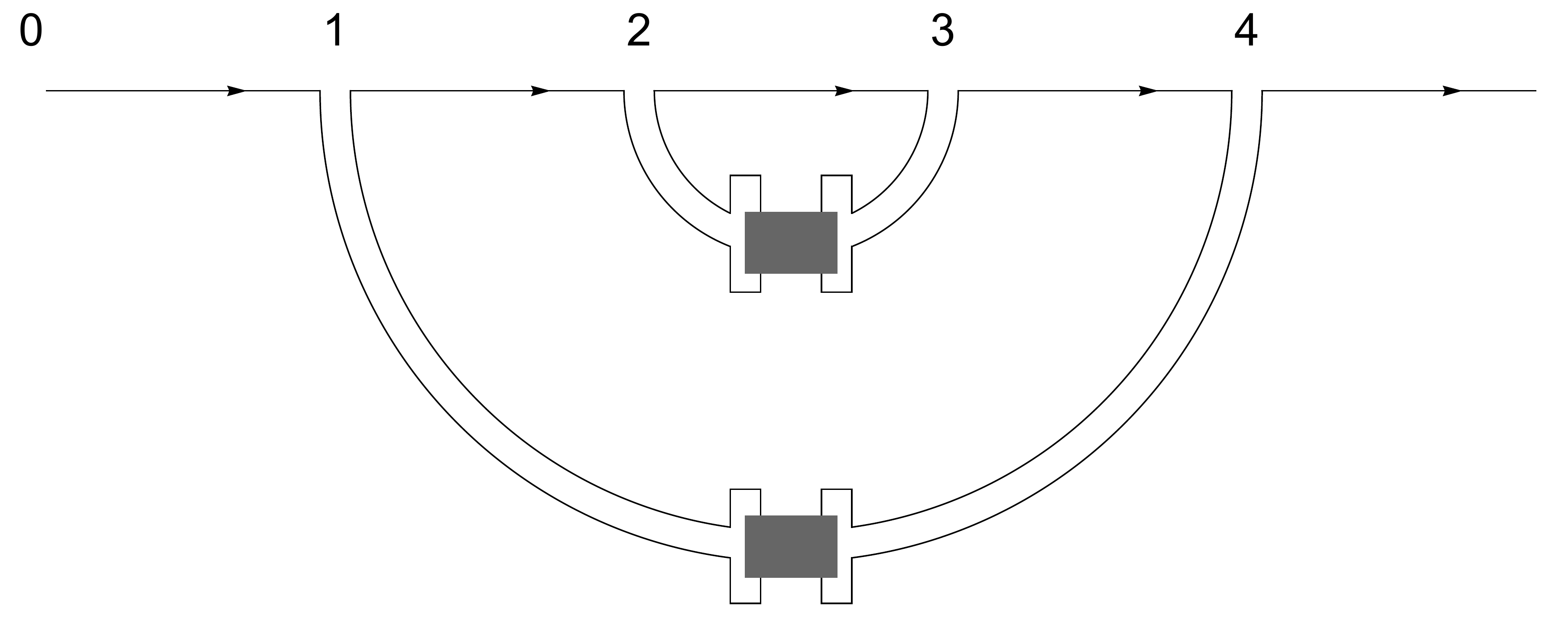}
\includegraphics[width=5.8cm]{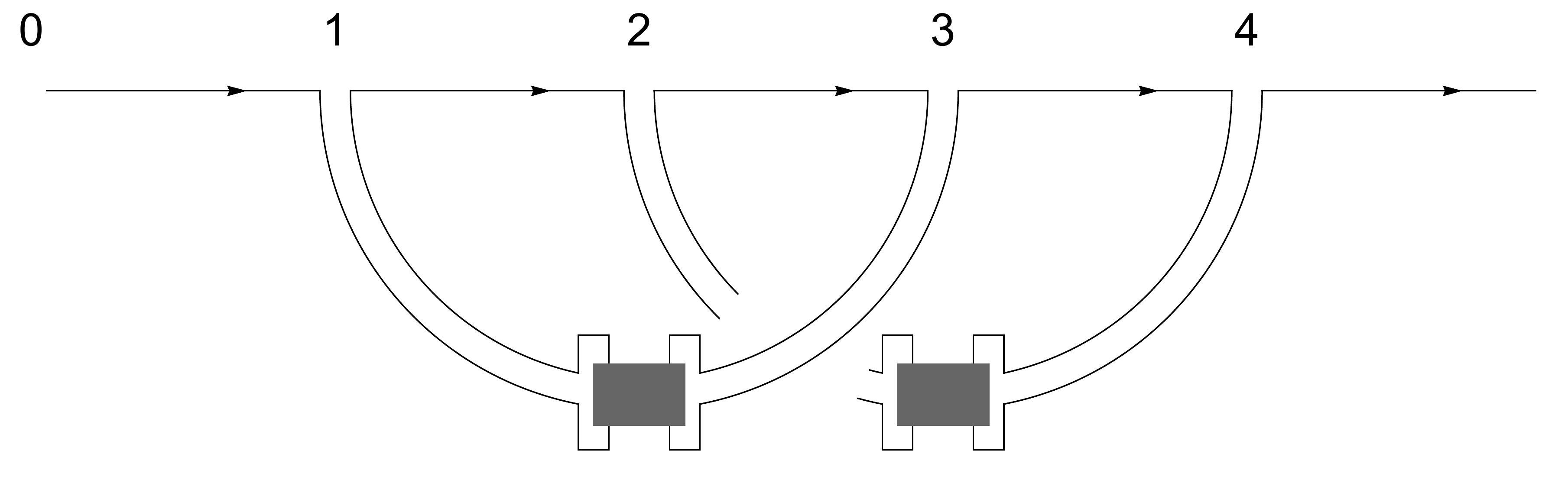}
\caption{The three basic configurations in circuit topology: Series, Parallel and~Cross.}
\label{fig.SPX}
\end{figure}

The starting point of circuit topology is to summarise the properties of folded chains in terms of these pairwise relationships of contacts. A~pair of contacts can interact in precisely three ways that are shown in Figure~\ref{fig.SPX} called
$\mathcal{S} = \{\{1,2\},\{3,4\}\}$ (series), \mbox{$\mathcal{P} = \{\{1,4\},\{2,3\}\}$} (parallel) and $\mathcal{X} = \{\{1,3\},\{2,4\}\}$ (cross). 
%\linenumbers

For more complicated chains, such as $\mathcal{F}$, we can record the relationship between each pair of contacts, ignoring the rest of the contacts. For~$\mathcal{F}$, we can describe the results in a table like the one given below. The~order of first appearance along the chain together with the pairwise information is sufficient to reconstruct the whole chord~diagram.

\noindent \setlength{\tabcolsep}{7.9mm}\begin{tabular}{cccccc}
\toprule
\textbf{\boldmath{$\mathcal{F}$}}&\textbf{\{1,4\}}&\textbf{\{2,7\}}&\textbf{\{3,5\}}&\textbf{\{6,9\}}&\textbf{\{8,10\}}\\
\midrule
\{1,4\}&0&X&X&S&S\\
\{2,7\}&X&0&P&X&S\\
\{3,5\}&X&P&0&S&S\\
\{6,9\}&S&X&S&0&X\\
\{8,10\}&S&S&S&X&0\\
\bottomrule
\end{tabular}

\vspace{+15pt}
The matrix of partial relationships is more concise and was shown to be a good way to describe the behaviour of folded chains in practice~\cite{Ma21,MWT14}.

Later in the paper, Section~\ref{sec.gCT}, we will use the above notation to introduce generalized circuit topology diagrams where each pair comes with a subscript $0,1$ or $-1$.
If the subscript is $0$ then the contact is hard as above, while the $\pm 1$ indicate a soft contact in the form of a left or right handed clasp.
For example the diagram that appears in Figure~\ref{fig.gCTExample1} is described as $\{\{1,4\}_1,\{2,7\}_{-1},\{3,5\}_0,\{6,9\}_0,\{8,10\}_{-1}\}$

The circuit topology description of folded chains can thus be generalised to also
include knotting and tangling effects, see~\cite{GM20}. We will see that any knotted folded chain can be described by a similar set of pairs of integers partitioning $1,2\dots 2n$ together with an index to specify the nature of the contact: hard or~soft.

\section{Knots and~Tangles}

Even when a chain does not bond to itself it can still be entangled, and provided we keep the endpoints fixed, such entanglement can constrain the chain as much as hard contacts would.
There is a large literature on knot theory that deals with this special case and in the next two sections we briefly review some aspects that are necessary for our purposes. 
The reader is warned that we use slightly non-standard versions of the usual definitions. We start with an example and a definition of tangle~diagrams.

\begin{figure}[H]
\includegraphics[width=7cm]{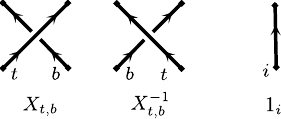}
\caption{From left to right: the positive crossing $X_{t,b}$, the~negative crossing $X^{-1}_{t,b}$ and the standard crossingless strand $1_i$. }
\label{fig.Xings}
\end{figure}

\begin{Definition}
A {\bf tangle diagram}
 is a finite directed graph whose edges are a disjoint union of oriented open paths called {\bf strands}. Each strand carries a distinct label. Strands are allowed to meet only in special vertices called crossings that look like like one of the models shown in Figure~\ref{fig.Xings}. 
\end{Definition}
\vspace{-6pt}

\begin{figure}[H]
\includegraphics[width=8cm]{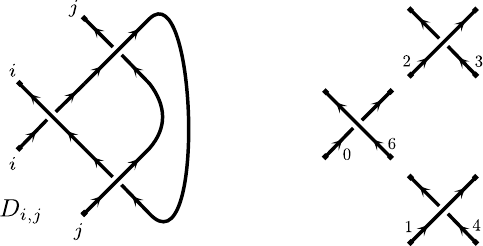}
\caption{Left: The tangle diagram $D_{i,j}$ with two strands called $i$ and $j$. $D$ can be assembled from the three crossings shown on the right by merging strands $0,2,4,6$ and $1,3$ as explained~below.}
\label{fig.TangleDiagramD}
\end{figure}

For a tangle diagram $D$, we often list the labels of its strands as a subscript. For~example the diagram $D_{i,j}$ shown in Figure~\ref{fig.TangleDiagramD} has two strands named $i$ and $j$. In~drawing pictures the labels are sometimes~suppressed.

The simplest tangle diagram is perhaps the crossing-less strand with label $i$, denoted by $1_i$.
Next the positive crossing $X_{t,b}$ with $t$ labeling the strand that is on top and $b$ marking the strand that is at the bottom. Likewise the negative crossing with top strand $t$ and bottom strand $b$ is denoted $X^{-1}_{t,b}$, see Figure~\ref{fig.Xings}.

\begin{figure}[H]
\includegraphics[width=12cm]{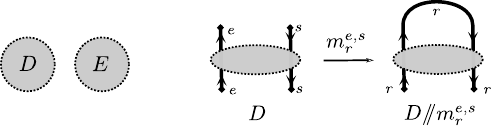}
\caption{\textbf{Left}: disjoint union of diagrams $D,E$. \textbf{Right}: merging strands $e,s$ of diagram $D$.}
\label{fig.Operations}
\end{figure}

Any more complicated tangle diagram can be assembled from the above examples using the fundamental operations on diagrams shown in Figure~\ref{fig.Operations}: disjoint union and strand merging. This provides a convenient notation for tangle diagrams allowing us to compute and discuss them~recursively.
\vspace{-6pt}

First we denote by $DE$ the {\bf disjoint union} of diagrams $D$ and $E$. Second, denote by $D\pp m^{e,s}_r$ the diagram obtained from $D$ by {\bf merging} the end of strand $e$ with the start of strand $s$, labelling the resulting strand $r$. By~merging two strands we mean connecting them without introducing new crossings. $D$ is also assumed to have strands labelled $e,s$ and no strand labelled $r$. 
Mergings may not always be defined but when the ends are close together merging them does make~sense.

For example the two strand tangle diagram $D_{i,j}$ shown in Figure~\ref{fig.TangleDiagramD} is denoted by
\begin{equation}
 D_{i,j} = X^{-1}_{6,0}X_{1,4}X_{2,3}\pp m^{0,2}_0\pp m^{0,4}_0\pp m^{0,6}_i\pp m^{1,3}_j.   
\end{equation}

\begin{figure}[H]
\includegraphics[width=6cm]{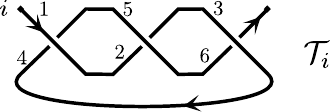}
\caption{The trefoil knot $\mathcal{T}_i$ viewed as a tangle with a single strand labelled $i$. The~numbers next to the crossings suggest how $\mathcal{T}_i$ can be merged from three~crossings.}
\label{fig.Trefoil}
\end{figure}

The special case where a tangle diagram has a single strand is called a knot diagram, see, for example, Figure~\ref{fig.Trefoil}. Often the ends of a knot diagram are required to be connected to form a loop but we prefer to leave them~open.

\begin{figure}[H]
\includegraphics[width=10cm]{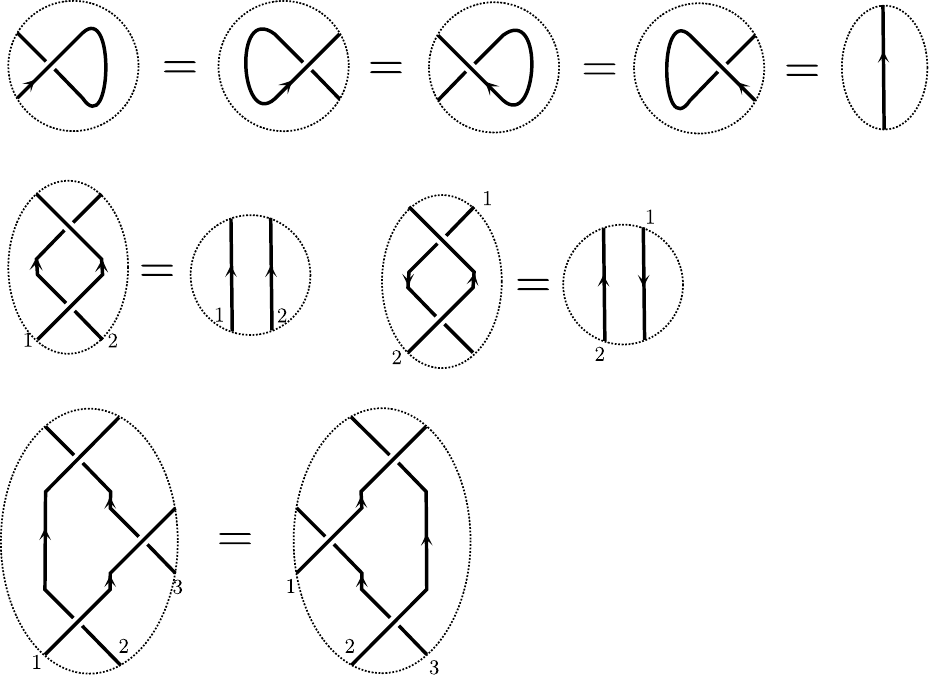}
\caption{The Reidemeister moves. Two tangles that locally differ as in one of the pictures are considered to be~equivalent.}
\label{fig.Reidemeister}
\end{figure}

A pair of tangle diagrams is said to be {\bf equivalent} if they can be made equal by repeatedly making local replacements as shown in Figure~\ref{fig.Reidemeister}. These replacements are commonly known as Reidemeister moves. The~theorem of Reidemeister~\cite{BZ03} asserts that tangle diagrams represent the same knot if and only if they are equivalent in the above sense. There are many more variants of the moves shown in the figure but these can all be realized as combinations of the ones shown~\cite{P10}.

Finally we introduce a short-hand notation to avoid long formulas
\begin{equation}
D\pp m^{a,b,c,d,\dots}_k = D\pp m^{a,b}_k \pp m^{k,c}_k \pp m^{k,d}_k \pp \dots .
\end{equation}

For example the diagram for the trefoil knot $\mathcal{T}_i$ shown in Figure~\ref{fig.Trefoil} can be described as
\begin{equation}
\mathcal{T}_i=X_{1,4}X_{5,2}X_{3,6}\pp m^{1,2,3,4,5,6}_i.
\end{equation}

\vspace{-12pt}

\section{Gamma-Calculus and the Alexander~Polynomial}

The fundamental difficulty of studying knots in terms of their diagrams is that there are many diagrams that represent the same knot. This is made precise in terms of the Reidemeister moves introduced in the previous section, see Figure~\ref{fig.Reidemeister}. To~obtain information about the knot using a diagram we need to compute a quantity that does not change when a Reidemeister move is applied. Such quantities are known as {\bf knot invariants}. 

In this section, we introduce one of the most useful knot invariants, called the {\bf Alexander polynomial}. There are many ways to describe and compute the Alexander polynomial. Here we choose a recursive definition in terms of tangle diagrams that fits well with our notation for tangle diagrams. It is known as $\Gamma$ calculus developed by Bar--Natan, see~\cite{BS13, Vo18} and~\cite{BN15} (Section~9) for more~context. 

\begin{Definition}[\textbf{\boldmath{$\Gamma$-calculus}}]
\label{def.Gamma}
For a tangle diagram $D$ whose strands are labelled by set $L$ define $\Gamma(D)=(\omega,A)$ 
where $A = \sum_{i,j\in L} A_{ij}r_ic_j$ and for each $i,j\in L$ the coefficients $A_{ij}$ and $\omega$ are rational functions in $t$.
\begin{enumerate}
    \item $\Gamma(1_i) = (1,r_ic_i)$.
    \item $\Gamma(X^{\pm}_{ij})=(1,r_ic_i+(1-t^{\pm 1})r_ic_j+t^{\pm 1}r_jc_j)$.
    \item For $(\omega,A)$ as above define
\begin{equation}
    (\omega,A)\pp \mu^{e,s}_k = \Big((1-A_{es})\omega,A+\frac{(\p_{r_e}A)(\p_{c_s}A)}{1-A_{es}}\Big)|_{\substack{r_e,c_s\mapsto 0\\ r_s,c_e\mapsto r_k,c_k}}
    \end{equation}
    Then $\Gamma(D\pp m^{e,s}_k) = \Gamma(D)\pp \mu^{e,s}_k$.
    \item If $\Gamma(D') = (\omega',A')$ then $\Gamma(D,D') = (\omega\omega',A+A')$.
\end{enumerate}
\end{Definition}

\begin{figure}[H]
\includegraphics[width=12cm]{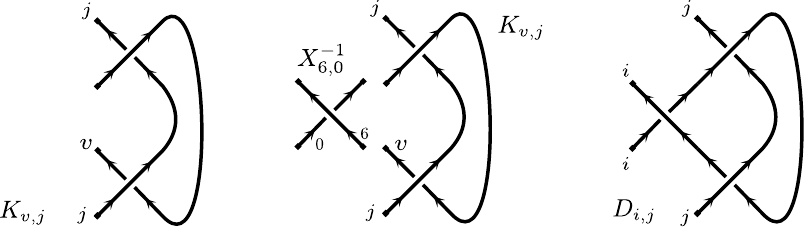}
\caption{Constructing the tangle $D$ from the tangle $K$ and a single~crossing.}
\label{fig.ConstructingD}
\end{figure}

Here, $\p_{r_e}A$ means the partial derivative of $A$ with respect to the variable $r_e$ and the notation $r_e,c_s\mapsto 0$ means we should set both $r_e$ and $c_s$ to $0$
and $r_s$ to $r_k$ and $c_e$ to $c_k$.
To see how this works in practice let us compute $\Gamma(D_{i,j})$ where $D_{ij}$ is the tangle shown in \mbox{Figure~\ref{fig.TangleDiagramD}}. 
It is efficient to start with the two-crossing tangle $K_{v,j} = X_{1,4}X_{2,3}\pp m^{1,3}_j\pp m^{2,4}_v$ and then
build $D_{ij}$ by attaching the third crossing to $K$ as follows (see also Figure~\ref{fig.ConstructingD})
\begin{equation}
 D_{i,j} = X^{-1}_{6,0}K_{v,j}\pp m^{0,v}_0\pp m^{0,6}_i.   
\end{equation}
\vspace{-12pt}

To compute $\Gamma(K_{v,j})$ we first compute $\Gamma(X_{1,4}X_{2,3})= $
\begin{equation*}
(1,r_1c_1+(1-t)r_1c_4+tr_4c_4+r_2c_2+(1-t)r_2c_3+tr_3c_3) = (\omega,A),
\end{equation*}
using parts (2) and (4) of Definition \ref{def.Gamma}. Then part (3) with $e,s = 1,3$ tells us we should compute $(\p_{r_1}A)(\p_{c_3}A) =(c_1+(1-t)c_4)((1-t)r_2+tr_3)$ and $A_{13} = 0$ so
\begin{equation*}
\Gamma(X_{1,4}X_{2,3}\pp m^{1,3}_j) = \big(1,A+(c_1+(1-t)c_4)((1-t)r_2+tr_3)\big)|_{\substack{r_1,c_3\mapsto 0\\ r_3,c_1\mapsto r_j,c_j}} = 
\end{equation*}
\begin{equation*}
(1,(1-t)r_2c_j+r_2c_2+(1-t)^2r_2c_4+tr_jc_j+t(1-t)r_jc_4+tr_4c_4) = (1,B).
\end{equation*}

Applying part (3) once more with $e,s = 2,4$ we find $B_{24} = (1-t)^2$ and $(\p_{r_2}B)(\p_{c_4}B) = ((1-t)c_j+c_2+(1-t)^2r_2)((1-t)^2r_2+t(1-t)r_j+tr_4)$ so
\begin{equation*}
\Gamma(K_{v,j}) = \Gamma(X_{1,4}X_{2,3}\pp m^{1,3}_j\pp m^{2,4}_v) = 
\end{equation*}
\[(2t-t^2, B+ \frac{((1-t)c_j+c_2+(1-t)^2c_4)((1-t)^2r_2+t(1-t)r_j+tr_4)}{2t-t^2})_{\substack{r_2,c_4\mapsto 0\\ r_4,c_2\mapsto r_v,c_v}}\]
\begin{equation*}
=(2t-t^2, tr_jc_j+ \frac{((1-t)c_j+c_v)(t(1-t)r_j+tr_v)}{2t-t^2}) 
\end{equation*}
\[
 = (2t-t^2,\frac{r_jc_j+(1-t)r_jc_v+(1-t)r_vc_j+r_vc_v}{2-t}).
\]

Now we are in a position to compute $\Gamma(D_{ij})$ in just three more steps, taking disjoint union with $X^{-1}_{6,0}$ and merging twice:
$\Gamma(X^{-1}_{6,0}K_{vj})=$
\begin{equation*}
(2t-t^2,r_6c_6+(1-t^{- 1})r_6c_0+t^{- 1}r_0c_0+\frac{r_jc_j+(1-t)r_jc_v+(1-t)r_vc_j+r_vc_v}{2-t}) = (\omega',A').
\end{equation*}

Since $A'_{0,v} = 0$ and $(\p_{r_0A'})(\p_{c_vA'}) = \frac{t^{- 1}c_0}{2-t}((1-t)r_j+r_v)$ we find \vspace{+6pt}
\end{paracol}
\nointerlineskip
\begin{equation*}
\Gamma(X^{-1}_{6,0}K_{vj}\pp m^{0,v}_0)=(2t-t^2,r_6c_6+(1-t^{- 1})r_6c_0+\frac{r_jc_j+(1-t)r_0c_j+t^{- 1}c_0((1-t)r_j+r_0)}{2-t}).
\end{equation*}
\begin{paracol}{2}
%\linenumbers
\switchcolumn

Finally, merging for the last time, the determined reader will find $\Gamma(D_{ij})=$  \vspace{+6pt}
\end{paracol}
\nointerlineskip
\begin{equation*}
\Gamma(X^{-1}_{6,0}K_{vj}\pp m^{0,v}_0\pp m^{0,6}_i)=(2t-t^2,(1-t^{- 1})r_ic_i+\frac{r_jc_j+t^{- 1}c_i(1-t)r_j+(1-t)c_jr_i+t^{-1}c_ir_i}{2-t}).
\end{equation*}
\begin{paracol}{2}
%\linenumbers
\switchcolumn

While instructive to do by hand once, these calculations tend to get a little tedious, so we recommend the reader to make use of a computer algebra package to do the calculations.
See the Appendix~\ref{sec.Appendix} for a description of the algorithm in Mathematica. As~the computations are conceptually very simple (only five lines of code!) it should be easy to implement the same algorithm in any other suitable~language. 
 
Another interesting example is that of the trefoil $\mathcal{T}_i = X_{1,4}X_{5,2}X_{3,6}\pp m^{1,2,3,4,5,6}_i$. In~precisely the same way as above we compute
\begin{equation}
\Gamma(\mathcal{T}_i)=(t^3-t^2+t,r_ic_i).
\end{equation}

We summarize the main properties of the $\Gamma$ calculus in the~following:

\begin{Theorem}[\textbf{\boldmath{$\Gamma$ is a tangle invariant generalizing Alexander \cite{BN15}}}]
Up to multiplication by a factor $\pm t^k$ in the $\omega$-part, $\Gamma$ is an invariant of tangles. For~tangles $K_i$ with one strand (i.e., knots)
we have $\Gamma(K_i) = (\Delta_K(t),r_ic_i)$, where $\Delta_K$ is the Alexander polynomial of knot $K$. 
\end{Theorem}

To understand how the factors $\pm t^k$ arise we note that the Alexander polynomial of knots has the same ambiguity. For~$\Gamma$ these factors arise because $\Gamma$ does not give the same value on both sides of a Reidemeister 1 move. For~example $\Gamma(1_i) = (1,r_ic_i)$ while $\Gamma(X_{12}\pp m^{12}_i)=(t,r_ic_i)$, see Figure~\ref{fig.Reidemeister1}.
\begin{figure}[H]
\includegraphics[width=6cm]{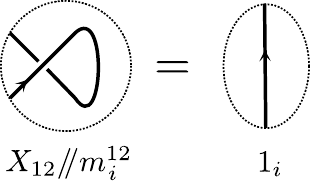}
\caption{Writing one of the Reidemeister 1 moves~algebraically.}
\label{fig.Reidemeister1}
\end{figure}

The other Reidemeister moves can likewise be written algebraically and it can be checked that $\Gamma$ yields the same result on both sides of each of the tangle equations representing Reidemeister 2 and 3.
For example, Reidemeister 3 is written as:
\begin{equation}
X_{12} X_{43} X_{56}\pp m^{14}_1\pp m^{25}_2\pp m^{36}_3 =  X_{16} X_{23} X_{45} \pp m^{14}_1\pp m^{25}_2\pp m^{36}_3,
\end{equation}
and it is easy to check (by computer) that indeed,
\begin{equation}
\Gamma(X_{12} X_{43} X_{56}\pp m^{14}_1\pp m^{25}_2\pp m^{36}_3) =  \Gamma(X_{16} X_{23} X_{45} \pp m^{14}_1\pp m^{25}_2\pp m^{36}_3).
\end{equation}

\section{\textbf{\boldmath{$H$}}-Tangle~Diagrams}
\label{sec.Htangles}

After warming up with knots we are now ready to introduce the topological framework for discussing folded chains: $H$-tangles.
These are an extension of the tangles from the previous sections where the strands may now bond to each other, forming what we call hard contacts 
or $H$-contacts. The~local structure of a $H$-contact is shown in Figure~\ref{fig.HContact} (Left).

\begin{figure}[H]
\includegraphics[width=12cm]{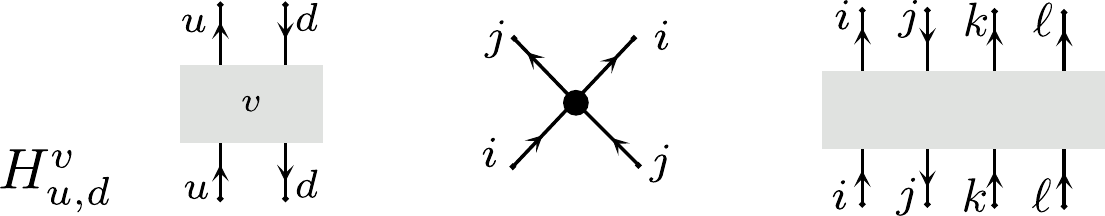}
\caption{\textbf{Left}: The vertex or hard contact in a $H$-tangle diagram with two strands labelled $u$ and $d$. \textbf{Middle}: a singular crossing in singular knots. \textbf{Right}: a generalized $H$-vertex.}
\label{fig.HContact}
\end{figure}

\begin{Definition}
A $H$-tangle diagram is the same as a tangle diagram except that one more type of vertex is allowed called a $H$-contact, see Figure~\ref{fig.HContact} (Left). Each $H$-contact carries a label.
Our notation for the $H$-vertex is $H^v_{u,d}$ where $u$ is the strand pointing upwards in the picture and
$d$ is the other strand. The~superscript $v$ denotes the label of the vertex.
\end{Definition}

One convenient way of labelling the $H$-contacts on a chain is by their first occurrence along the strand. When the $H$-contacts are not explicitly labeled, this labeling is assumed.
An example of a $H$-tangle with one strand already appeared in Figure~\ref{fig.gCTExample1} in the introduction. Another example with one strand and two hard contacts appear in Figure~\ref{fig.E}.
\begin{figure}[H]
\includegraphics[width=5cm]{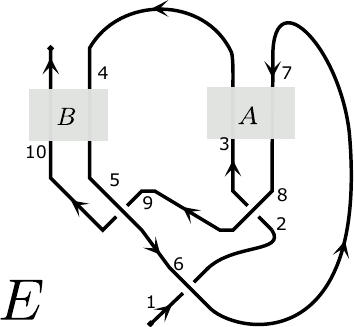}
\caption{The $H$-tangle diagram $E$ with two $H$-contacts labelled $A,B$.}
\label{fig.E}
\end{figure}

The notion of $H$-tangle diagrams is very similar to that of singular knots and tangles~\cite{BEHY13} but not exactly the same.
The difference is in how the strands pass through the vertex as is seen by comparing Figure~\ref{fig.HContact} Left ($H$-vertex) and Middle (singular vertex).
We also note that more complicated vertices can be modelled similarly to the proposed $H$-vertex, see for example, Figure~\ref{fig.HContact} Right, but we will leave
that for future~work.

Most of the theory for the usual tangle diagrams as developed in the previous sections is still valid. For~example, $H$-tangle diagrams can be assembled from the basic ones using disjoint union and merging. These operations are defined in precisely the same way as in the tangle~case.

\textls[+15]{For example, coming back to our discussion of circuit topology the basic configurations} $\mathcal{S},\mathcal{P},\mathcal{X}$ from Figure~\ref{fig.SPX} can now be interpreted as $H$-tangles as follows: \linebreak
\mbox{$\mathcal{S}_i=H_{1,2}H_{3,4}\pp m^{1,3,4,2}_i$} and
$\mathcal{P}_i=H_{1,2}H_{3,4}\pp m^{1,2,3,4}_i$.
The cross $\mathcal{X}_i$ is a little more involved as it involves crossings and there are in fact several options, one of which is:
\begin{equation}
\mathcal{X}_i = H_{12} H_{34}X_{12,5}X^{-1}_{14,7}X_{16,9}X^{-1}_{10,11}\pp m^{1,5,7,3,9,11,10,12,2,14,16,4}_i.    
\end{equation}

The example $E$ in Figure~\ref{fig.E} is likewise expressed as:
\begin{equation}
E_i = H^A_{10,4} H^B_{3,7}X_{6,1}X^{-1}_{8,2}X^{-1}_{5,9}\pp m^{1,2,3,4,5,6,7,8,9,10}_i.    
\end{equation}

The most important part of our model of $H$-tangles is the type of Reidemeister move we allow.
We model the hard contact $H^v_{u,d}$ as a rigid square and so the only way it is allowed to interact with other strands is for strands to pass over or under the vertex.
Two vertices are assumed to never be precisely on top of each other, just like crossings in a usual tangle diagram are not supposed to coincide.
More formally, we say two $H$-tangle diagrams are {\bf equivalent} if they can be made equal by repeatedly making local replacements as shown in Figure~\ref{fig.HReidemeister}.
We name these the Reidemeister 4~moves. 

Since the $H$-contacts are assumed to be rigid we should remember to always bring the $H$-contacts in the correct position with one up and one down, rotating in the counter-clockwise direction. Contacts with other orientations can be modelled using a $H$-contact by adding crossings as shown in Figure~\ref{fig.HXings}.
\begin{figure}[H]
\includegraphics[width=10cm]{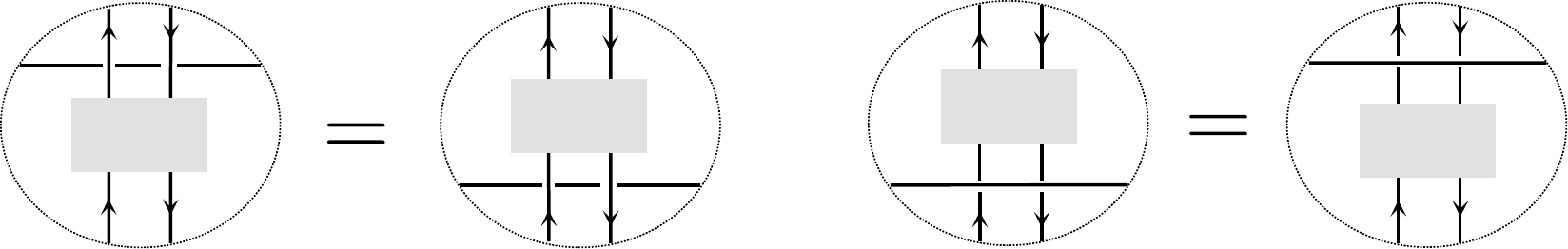}
\caption{New Reidemeister 4 moves for $H$-tangle diagrams. The~orientation of the horizontal strand is arbitrary so in fact we have four such~moves.}
\label{fig.HReidemeister}
\end{figure}

\vspace{-6pt}

\begin{figure}[H]
\includegraphics[width=6cm]{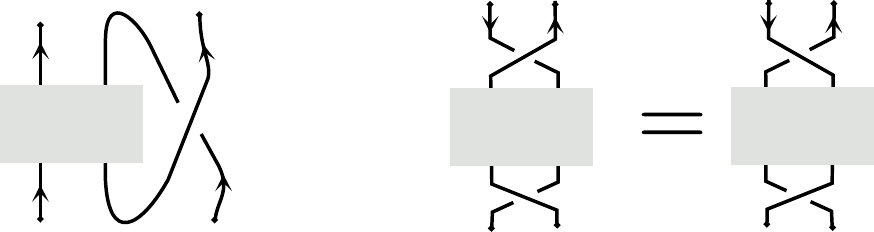}
\caption{In case we want two strands to run in parallel near a $H$-contact we augment it with a crossing. We can also flip it over using either of the two pictures on the~right. }
\label{fig.HXings}
\end{figure}

It may not be immediately clear that the two pictures on the right of Figure~\ref{fig.HXings} represent the same $H$-tangle. This is demonstrated in the next Figure~\ref{fig.HWhitney}.
The first step uses Reidemeisters 1--3 while the second uses Reidemeister 4 to pass the horizontal strands under the $H$-contact. The~next equality is obtained by sliding the $H$-contact along the loop
and the final step is like the~first.
\begin{figure}[H]
\includegraphics[width=9cm]{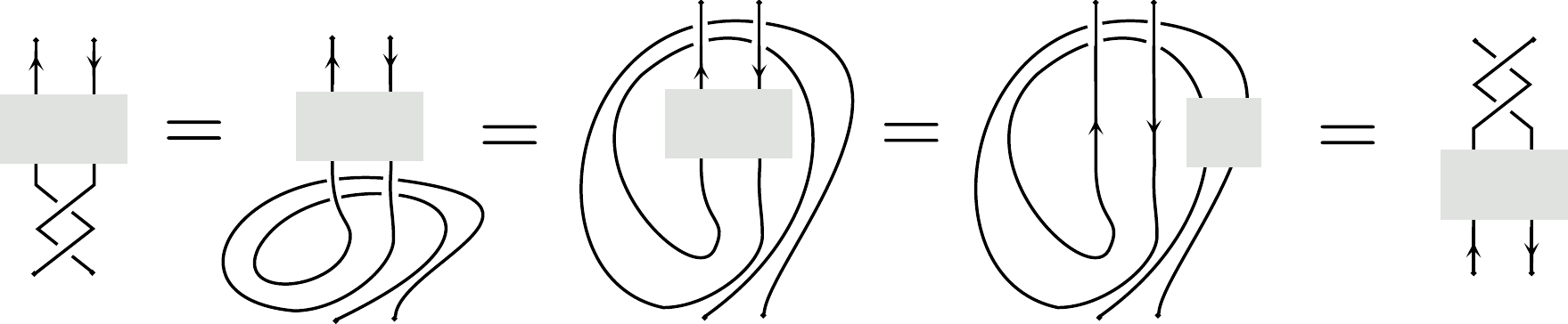}
\caption{A proof that two crossings can be passed through an $H$-contact using Reidemeisters 1--4.}
\label{fig.HWhitney}
\end{figure}

It should be noted that many variations on the notion of $H$-tangle introduced here are possible by adding more Reidemeister moves for the $H$-contact to satisfy.
Each additional Reidemeister move would make the vertex more flexible. For~example one could assume $H_{ij} = H_{ji}$ which would give a theory equivalent to that of
singular tangles. Indeed, this assumption would correspond to the fifth move for singular knots, see $\Omega5a$ on p. 3 of~\cite{BEHY13}. 
Another natural symmetry one could require is that the left-hand diagram of Figure~\ref{fig.HXings} should be equivalent to that where the crossing is opposite.
Since it is conceptually simplest we prefer to not assume any Reidemeister moves beyond the 4-th. In~principle it is not difficult to adjust the theory by including an 
additional Reidemeister move but we leave this for future~work.

\section{The Double Gamma Invariant of \textbf{\boldmath{$H$}}-Tangles}

In this section, we aim to extend $\Gamma$-calculus to $H$-tangles. We first try the most direct approach and show its limitations. Next we
propose a more powerful version of $\Gamma$ called $\DG$ that is more suitable for working with $H$-tangles.
 
The key idea of $\Gamma$ calculus is to assign a pair $(\omega,A)$ to any tangle where $A$ has as many $r$-variables as the number strands in the tangle.
Since $\Gamma(H_{ij})$ involves two strands, we assign to it a value of the form:
\begin{equation}
\Gamma(H^a_{ij}) = 
(p_0,p_1 r_ic_i+p_2 r_ic_j+p_3 r_jc_i+p_4r_jc_j),    
\end{equation}
for some functions $p_0,p_1,p_2,p_3,p_4$ of $t$ and $a$ that are yet to be~determined.

To ensure that $\Gamma$ takes the same value on equivalent $H$-tangle diagrams, it suffices
to choose the value $\Gamma(H^a_{ij})$ so that the Reidemeister 4 equations are satisfied.
The Reidemeister 4 moves appeared in Figure~\ref{fig.HReidemeister}, here we wrote them algebraically and applied $\Gamma$ to both~sides:

\begin{align*}
\Gamma(H^a_{12}X_{34}X^{-1}_{56}\pp m^{35}_3\pp m^{14}_1\pp m^{62}_2) &= \Gamma(H^a_{12}X_{34})X^{-1}_{56}\pp m^{35}_3\pp m^{41}_1\pp m^{26}_2)\\
\Gamma(H^a_{12}X^{-1}_{43}X_{65}\pp m^{35}_3\pp m^{14}_1\pp m^{62}_2) &= \Gamma(H^a_{12}X^{-1}_{43}X_{65}\pp m^{35}_3\pp m^{41}_1\pp m^{26}_2).
\end{align*}

Since we can explicitly compute both sides using $\Gamma$ calculus, it is left to the reader to solve these equations for $p_0,\dots p_4$
and verify (By computer would be most convenient!) that they are equivalent to $p_2=1-p_1, p_3=1-p_1, p_4=p_1$.

As a consequence we can compute what a single $H$-vertex connected to itself, see Figure~\ref{fig.HKink} (Left) evaluates to: 
$\Gamma(H^a_{1,2}\pp m^{12}_i) = (p_0p_1,r_ic_i)$. This is problematic because it follows that for any diagram $Y_j$
$\Gamma$ takes the same value on the two configurations in the middle and right of Figure~\ref{fig.HKink}. 
 In particular this implies that $\Gamma$ cannot distinguish the two $H$-tangles $\mathcal{S,P}$ representing two of the basic pair configurations in circuit topology (see Figure~\ref{fig.SPX}): 
$\mathcal{S}_1=H_{1,2}H_{3,4}\pp m^{1,3,4,2}_i$ and
$\mathcal{P}_1=H_{1,2}H_{3,4}\pp m^{1,2,3,4}_i$.
The same argument shows that the $\Gamma$ invariant is insensitive to the location in the chain of any element $\mathcal{B}_i$ with $\Gamma(\mathcal{B}_i)=(\omega,r_ic_i)$
\begin{figure}[H]
\includegraphics[width=7cm]{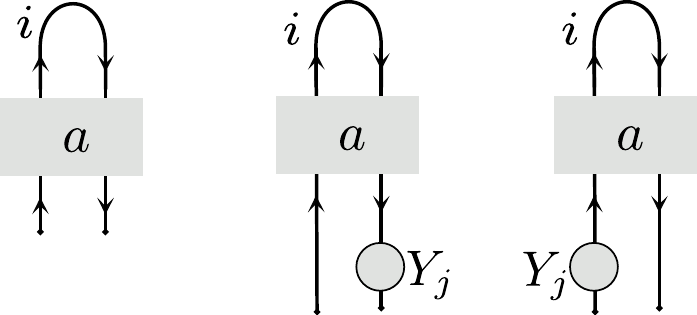}
\caption{\textbf{Left}: A $H$-vertex on a single strand cannot be detected by the naive extension of $\Gamma$ calculus. It takes the same value on both pictures on the~\textbf{Right}.}
\label{fig.HKink}
\end{figure}

We conclude that introducing the above $\Gamma$ invariant for $H$-tangles misses too many fundamental features of $H$-tangles such as the basic pair interactions.
This brings us to propose a more sophisticated version of $\Gamma$ called $\DG$ (double $\Gamma$) that distinguishes a lot more about $H$-tangles, while retaining a close relation to $\Gamma$.

\begin{Definition}[\textbf{\boldmath{$\DG$ calculus}}]
For a $H$-tangle diagram $D$ whose strands are labelled by a set $L$ of positive integers 
and $H$-vertices indexed by set $V$ define $\DG(D)$ by the following rules.
\mbox{$\DG(D)=(\omega,A)$} where $A = \sum_{\pm i,\pm j\in L} A_{ij}r_ic_j$ and $\omega,A_{ij}$ rational functions of $t,h_v,v\in V$.
\begin{enumerate}
    \item $\DG(1_i) = (1,r_ic_i+r_{-i}c_{-i})$.
    \item $\DG(X_{i,j}) = \big(1,r_ic_i+r_{-i}c_{-i}+t^2(r_jc_j+r_{-j}c_{-j})+(1-t)(tr_i+r_{-i})(c_{j}+c_{-j})\big)$.
\item $\DG(X^{-1}_{i,j}) = \big(1,r_ic_i+r_{-i}c_{-i}+t^{-2}(r_jc_j+r_{-j}c_{-j})+(1-t^{-1})(r_i+t^{-1}r_{-i})(c_{j}+c_{-j}))\big)$.
\item $\DG(H^B_{i,j}) = \big(1,t(r_{-i}c_{-i}+r_{j}c_{-j})+(1-t)(r_{-j}c_{-i}+r_{-i}c_{-j})+(1-h_B)(tr_{i}c_{-j}-tr_{-j}c_{-j}+r_{-j}c_{j})+r_{-j}c_{i}+h_Br_{i}c_{j}\big)$.
    \item $\DG( D\pp m^{i,j}_k) = (\omega,A)\pp \mu^{i,j}_k\pp \mu^{-i,-j}_{-k}$, where $\mu$ was defined in Definition \ref{def.Gamma}.
    \item If $\DG(D') = (\omega',A')$ then $\DG(D,D') = (\omega\omega',A+A')$.
\end{enumerate}
\end{Definition}

As the rules for $\DG$ are very similar to those of $\Gamma$-calculus, we do not work out any examples step by step but instead refer the reader to the computer code in the Appendix~\ref{sec.Appendix}.
For example for the $H$-tangle $E$ from Figure~\ref{fig.E} the value of $\DG$ is:
\[\DG(E_1) = (-h_A h_B+t h_A+h_B-t, (1-t)r_{-1}c_{-1} +r_{-1}c_1+tr_1c_{-1})\].

The problematic $H$-tangle $\mathcal{B}$ discussed at the start of this section that $\Gamma$ fails to handle properly 
also has a non-trivial $\DG$ value:
\[
\DG(\mathcal{B}_1) =(t(1-h_1), (1-t)r_{-1}c_{-1} +r_{-1}c_1+tr_1c_{-1}).
\]

The second part of $\DG(D) = (\omega,A)$ is often quite lengthly so we often just list the $\omega$-part and use the notation $\DG_1(D) = \omega$. In~fact this notation was already used in the introduction where $\DG_1$ of the $H$-tangle from Figure~\ref{fig.gCTExample1} was listed. As~expected $\DG_1$ also distinguishes the three basic pair interactions $\mathcal{S,P,X}$ of circuit topology. 
Labelling the $H$-contacts in order of appearance as one walks the chain we find
\begin{align*}
\DG_1(\mathcal{S}) &= (1-h_1)(1-h_2)\\
\DG_1(\mathcal{P}) &= \left(h_2-1\right) \left(-h_1 t+t-1\right)\\
\DG_1(\mathcal{X}) &= -h_1 \left(h_2-1\right) \left(t^3-t^2+1\right) t^2\\&+h_2 \left(t^4-t^3+2 t-1\right) t-t^5+2 t^4-3 t^3+3 t^2-3 t+2.
\end{align*}

As with the usual $\Gamma$ calculus, there is a minor ambiguity in $\DG_1$ due to fact that it is not invariant under Reidemeister 1. 
$\DG_1$ should always be taken up to multiplication by a factor $\pm t^k$. Actually the situation with Reidemeister 1 is a bit more serious
this time and to make $\DG$ an invariant we should always use {\bf 0-framed}
$H$-tangle diagrams. We say an $H$-tangle diagram is $0$-framed if after cutting out the $H$-vertices, all the remaining pieces of strand have writhe $0$. Recall that the writhe of a strand is the sum of the signs of all the crossings of the strand with itself. This is not a serious restriction as any diagram can be modified to have writhe $0$ by inserting curls where necessary using the Reidemeister 1~move.

\begin{Theorem}[\textbf{\boldmath{$\DG$ is an invariant of $H$-tangles}}]
If two $0$-framed $H$-tangle diagrams $D,E$ are equivalent then $\DG(D) = \DG(E)$ (up to a factor $\pm t^k$ in the $\DG_1$ parts). If $K$ is a knot diagram, then $\DG_1(K) = \Delta_K(t^2)$, that is,~the Alexander polynomial of the knot at $t^2$.
\end{Theorem}
\begin{proof}
Like ordinary $\Gamma$ calculus, the proof of invariance under the Reidemeister moves can be carried out by computer. One simply has to verify that
both sides of all Reidemeister moves yield the same result after applying $\DG$. We already sketched above that in the case of Reidemeister 1, one needs to correct for the framing and an ambiguity of a factor $\pm t^k$ will arise in the $\omega$ part. 

When no $H$-contacts are present it can be verified that $\DG$ is just $\Gamma$ applied to the double of a knot. By~the double, we mean that a parallel strand is added next to the existing strand. 
It is well-known that the Alexander polynomial of the doubled knot is the Alexander polynomial of the original taken at $t^2$, see for example Chapter 8 of~\cite{BZ03}.
\end{proof}

The invariant $\DG$ is new but as we said it is closely related to $\Gamma$-calculus, not just for knots. To~obtain the rules for $\DG$ calculus we applied $\Gamma$ calculus to the $H$-tangle in which each strand was doubled. The~value for the $H$-contact is however not of this type but rather found by the method of undetermined coefficients, imposing all the necessary Reidemeister moves. 
Imposing the Reidemeister moves does not determine all coefficients so the remaining ones were chosen in a way to at least distinguish the basic examples listed in the tables but not involve too many additional~variables.

\section{Generalized Ciruit~Topology}
\label{sec.gCT}

So far, we have have modelled folded chains using $H$-tangles that have a single strand and we developed a topological invariant $\DG$ of these objects.
In circuit topology a chain is described by its contacts. The~hard contacts are obvious in a $H$-tangle but the concept of soft contacts is less clear. 
In~\cite{GM20} four types of soft contacts were proposed as the four simplest ways a single strand can knot itself, see Figure~\ref{fig.SoftContacts}.
These coincide with all the possible knots one can tie in which the projection has only four~crossings.
\begin{figure}[H]
\includegraphics[width=11cm]{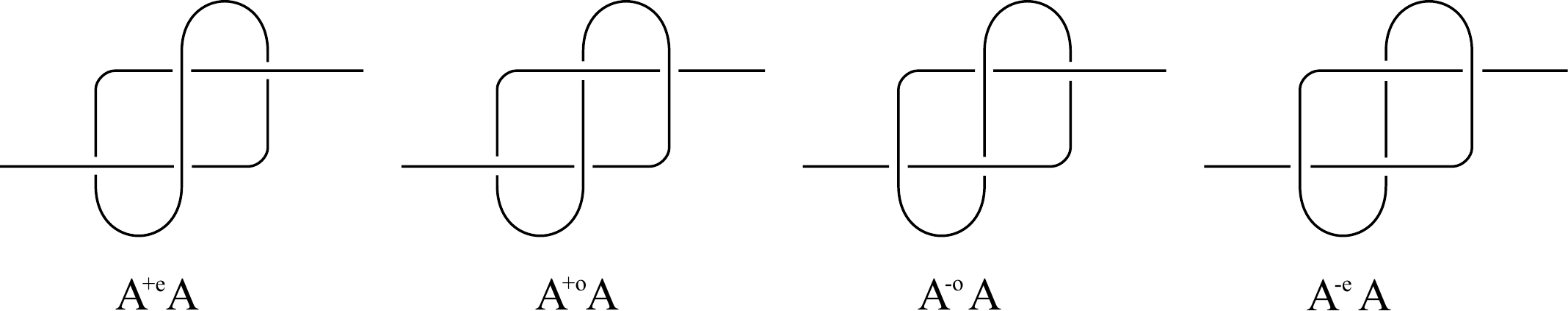}
\caption{The four soft contacts, reproduced from Figure~8 of~\cite{GM20}.}
\label{fig.SoftContacts}
\end{figure}

In this work, we approach the notion of soft contact from a tangle point of view. This means that we temporarily allow ourselves to think of the chain as
being composed of several strands. This simplifies matters because two strands can interlock non-trivially already with only two crossings. 
The simplest way two strands can interlock is the clasp, shown in Figure~\ref{fig.Clasps}. Like the crossing there is a positive and a negative version
of the clasp. In~what follows we will refer to either of them as a (soft) clasp~contact.  
\begin{figure}[H]
\includegraphics[width=7cm]{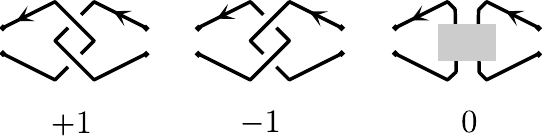}
\caption{The circuit topology contacts: \textbf{Left}: the positive and negative clasp, \textbf{Right}: the hard~contact.}
\label{fig.Clasps}
\end{figure}

The clasp contacts are universal in the sense that attaching a clasp next to a crossing causes the crossing to change sign, see Figure~\ref{fig.CrossingChange}. This means any other soft contact or really any tangle must be expressible using clasps only. For~example the four model configurations of soft contacts from~\cite{GM20}, Figure~\ref{fig.SoftContacts}, correspond precisely to the non-trivial combinations of two clasps. In~Figure~\ref{fig.SoftToClasp}, we sketch more precisely how the first of the four soft contacts may explicitly be converted to the standard clasp diagram form of Defintion \ref{def.gCT} given below. The~three others can be done in precisely the same~way.
\begin{figure}[H]
\includegraphics[width=8cm]{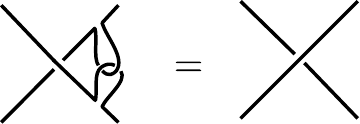}
\caption{Inserting a soft contact (clasp) has the effect of changing the sign of a~crossing.}
\label{fig.CrossingChange}
\end{figure}

\vspace{-6pt}
\begin{figure}[H]
\includegraphics[width=13cm]{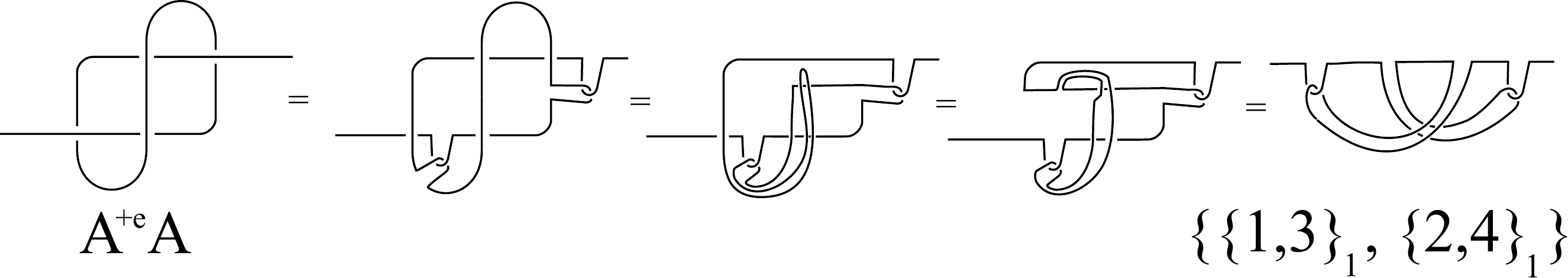}
\caption{Changing crossings by inserting clasps, straightening the picture keeping the clasps factors the soft contact $A^{+e}A$ into a gCT diagram with two~clasps.}
\label{fig.SoftToClasp}
\end{figure}

Even when restricting soft contacts to just be clasps the number of ways to combine clasps to get the same tangle is enormous making it hard to discern any
clear structure from them directly. To~talk about circuit topology more concretely we prefer a more rigid type of diagrams where the contacts interact
in a very simple way. We call such diagrams generalized circuit topology diagrams (gCT diagrams). In~the case of tangles, such diagrams appeared in
~\cite{MP19} under the name descending clasp diagrams. We basically write the chain as a horizontal line that is interrupted by arcs representing~contacts.

\begin{Definition}[\textbf{Generalized circuit topology (gCT) diagrams}]
\label{def.gCT}
A {\bf gCT diagram} is encoded by a partition 
\[\{\{a_1,b_1\}_{\sigma_1},\{a_2,b_2\}_{\sigma_2},\dots,\{a_n,b_n\}_{\sigma_n}\}\]
where $a_i,b_i$ are integers such that $\bigcup_{i=1}^n\{a_i,b_i\} = \{1,2,\dots 2n\}$ and $\sigma_i\in\{-1,0,1\}$.
Given this data define a one-strand $H$-tangle diagram
as follows. Starting with the $x$-axis oriented left to right, for~each pair $\{a,b\}_{\sigma}$ connect the points $a$ and $b$ with a contact of type $\sigma$ as shown in Figure~\ref{fig.Clasps2}. 
Whenever $a<c<b<d$ the lines corresponding to pairs $\{a,b\}$ and $\{c,d\}$ will cross the $\{a,b\}$ contact is assumed to be in front. No other crossings appear.
\end{Definition}
\vspace{-6pt}
\begin{figure}[H]
\includegraphics[width=11cm]{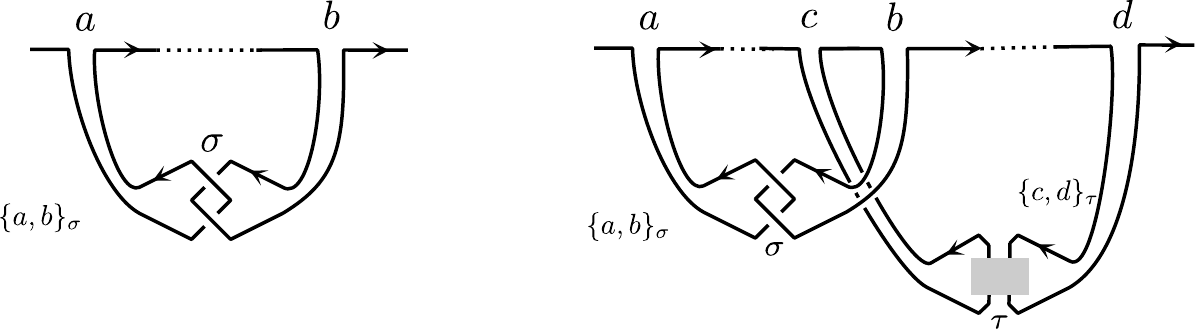}
\caption{The notation $\{a,b\}_\sigma$ means a contact of type $\sigma$ connecting $a,b$ below the $x$-axis. Two such contacts can only cross as~shown.}
\label{fig.Clasps2}
\end{figure}

Examples of gCT diagrams already appeared in Figures~\ref{fig.gCTExample1} and \ref{fig.CTExample}. The~notation for these diagrams is 
\[\{\{1,4\}_1,\{2,7\}_{-1},\{3,5\}_0,\{6,9\}_0,\{8,10\}_{-1}\}\]\[\{\{1,4\}_0,\{2,7\}_{0},\{3,5\}_0,\{6,9\}_0,\{8,10\}_{0}\}\] respectively.
Also the soft contacts from Figure~\ref{fig.SoftContacts} were factored into two claps using the reasoning in Figure~\ref{fig.SoftToClasp} and
in the gCT notation we thus find: 
\begin{align*}A^{+e}A &= \{\{1,3\}_{1},\{2,4\}_1\}\\
A^{+o}A &= \{\{1,3\}_{1},\{2,4\}_{-1}\}\\
A^{-o}A &= \{\{1,3\}_{-1},\{2,4\}_{1}\}\\
A^{-e}A &= \{\{1,3\}_{-1},\{2,4\}_{-1}\}.
\end{align*}

To get a sense of what gCT diagrams look like in general we show a randomly chosen gCT diagrams with 20 contacts in Figure~\ref{fig.gCTs}. In~gCT notation it is \vspace{+6pt}
\end{paracol}
\nointerlineskip
\[
\left\{\{1,14\}_1,\{2,11\}_0,\{3,13\}_0,\{4,12\}_1,\{5,8\}_{-1},\{6,9\}_1,\{7,19\}_{-1},\{10,18\}_0,\{15,20\}_{-1},\{16,17\}_0\right\}.
\]
\begin{paracol}{2}
%\linenumbers
\switchcolumn

\vspace{-6pt}
\begin{figure}[H]
\includegraphics[width=13.8cm]{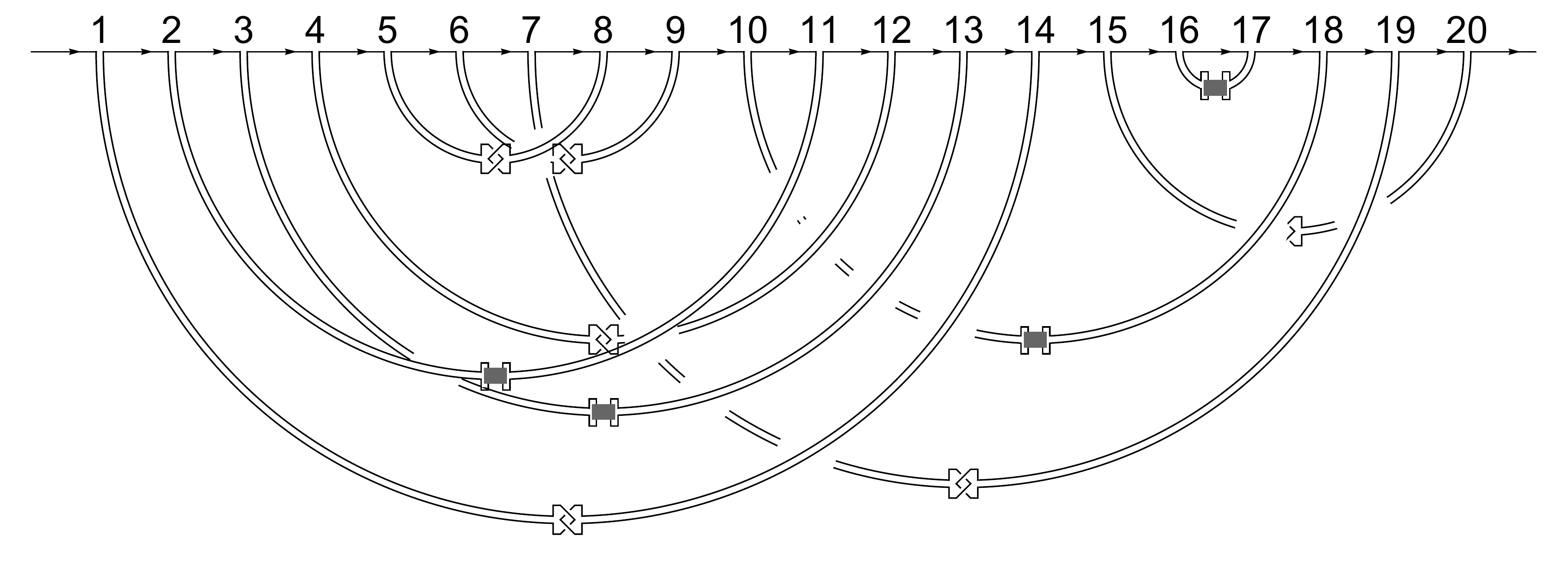}
\caption{A random gCT diagram with twenty~contacts. }
\label{fig.gCTs}
\end{figure}

As a justification for restricting to gCT diagrams, we claim that planar $H$-tangle diagrams can be represented this way.
By planar we mean that viewed as a decorated graph whose vertices are crossings and $H$-vertices it can be drawn in the plane.
All $H$-tangle diagrams
that are actual projections of three-dimensional chains are planar. Planar in the sense of the underlying graph, it can of course have crossings. 
Usually the planarity condition is part of the definition of a tangle diagram,
however since it is not needed for $\DG$ we preferred to work a little more~generally.

\begin{Theorem}[\cite{MP19}]
\label{thm.MP}
\textls[-15]{Any planar $H$-tangle diagram with a single strand is equivalent to a gCT~diagram.}
\end{Theorem}
This theorem is due to~\cite{MP19} although their argument is given for usual tangles they make it clear the same technique works for singular knots and also $H$-tangles.
A first hint of why it works is that inserting a clasp contact near a crossing changes the sign of that crossing as shown in Figure~\ref{fig.CrossingChange}.
This implies that we can freely change the crossings of our $H$-tangle diagram at the cost of inserting clasp contacts. This way all the complicated knots and tangles can be
'factored' into simple clasp contacts. The~hard contacts will remain and the clasp contacts may interlace in complicated ways but by changing some more crossings
we can simplify further until a gCT diagram remains. An~illustration of this argument in a very simple case already appeared in Figure~\ref{fig.SoftToClasp}.
A precise argument using the pure braid group would take us too far afield and for that and more we refer to~\cite{MP19}.

While convenient and interesting to study, the reader is warned that such diagrams are not necessarily unique. 
In fact in~\cite{MP19} a complete set of moves relating equivalent gCT diagrams is suggested in the knot case and we suspect this too holds in the $H$-tangle setting. Like the Reidemeister moves
themselves these moves are sufficiently complicated that it is still essential to compute knot invariants such as $\DG$ for distinguishing gCT~diagrams.

Many natural operations in both circuit topology and tangle theory can easily be handled in terms of gCT diagrams. To~give one simple example,
the $H$-$H$ interactions studied in~\cite{CEM21} can be described precisely as follows. Suppose we have two hard contacts that are interacting by a soft (clasp) contact as
Figure~\ref{fig.HH} Left. Assuming the two contacts are in series we can just slide the soft (clasp) contact off the arcs and turn it into an additional arc in the gCT diagram, shown on the right of the figure.
This means inserting soft (clasp) contacts like this can be done entirely on the level of gCT~diagrams.
\begin{figure}[H]
\includegraphics[width=12cm]{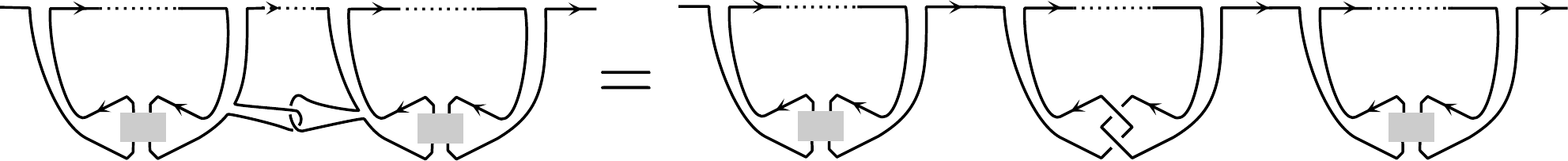}
\caption{Implementing a $H$-$H$ interaction by adding an additional soft (clasp) contact.}
\label{fig.HH}
\end{figure}

We end this section by listing all 27 gCT diagrams with at most two contacts, see Figure~\ref{fig.TablegCT2}. For~each diagram we also list its gCT notation.
The hard contacts of a gCT diagram are always labelled in order of first appearance as one travels along the strand. For~the $i$-th hard contact we used the variable $h_i$ in our expression
for $\DG_1$. Except~for the pair of mirror image trefoil knots $\{\{1,3\}_{1},\{2,4\}_{1}\}$ and $\{\{1,3\}_{-1},\{2,4\}_{-1}\}$ diagrams with the same value of $\DG_1$ are in fact equivalent.
It is well known that the Alexander polynomial does not distinguish mirror images of knots. The~ambiguity in multiplication by a power of $t$ has been used to make sure the lowest power of $t$ is always $t^0$.

The first eight diagrams in the table are all gCT diagrams on which $\DG_1$ takes the value $1$. Using Reidemeister moves or perhaps with the help of a rope the reader should be able to verify that each of them is in fact equivalent to the standard unknot diagram: a straight line with no contacts. Of~the diagrams that come after we now know for sure that they are not equivalent to the unknot as $\DG_1$ takes a value different from $1$.
 
In the Appendix~\ref{sec.Appendix} a similar table for all gCT diagrams with three contacts is given. These tables show that $\DG_1$ is able to distinguish a lot of relevant cases.
A table for all gCT diagrams with four contacts is also feasible but more and more equivalent gCT diagrams will appear and additional restrictions should be enforced to obtain a useful~table.

\begin{figure}[H]
\includegraphics[width=0.99\linewidth]{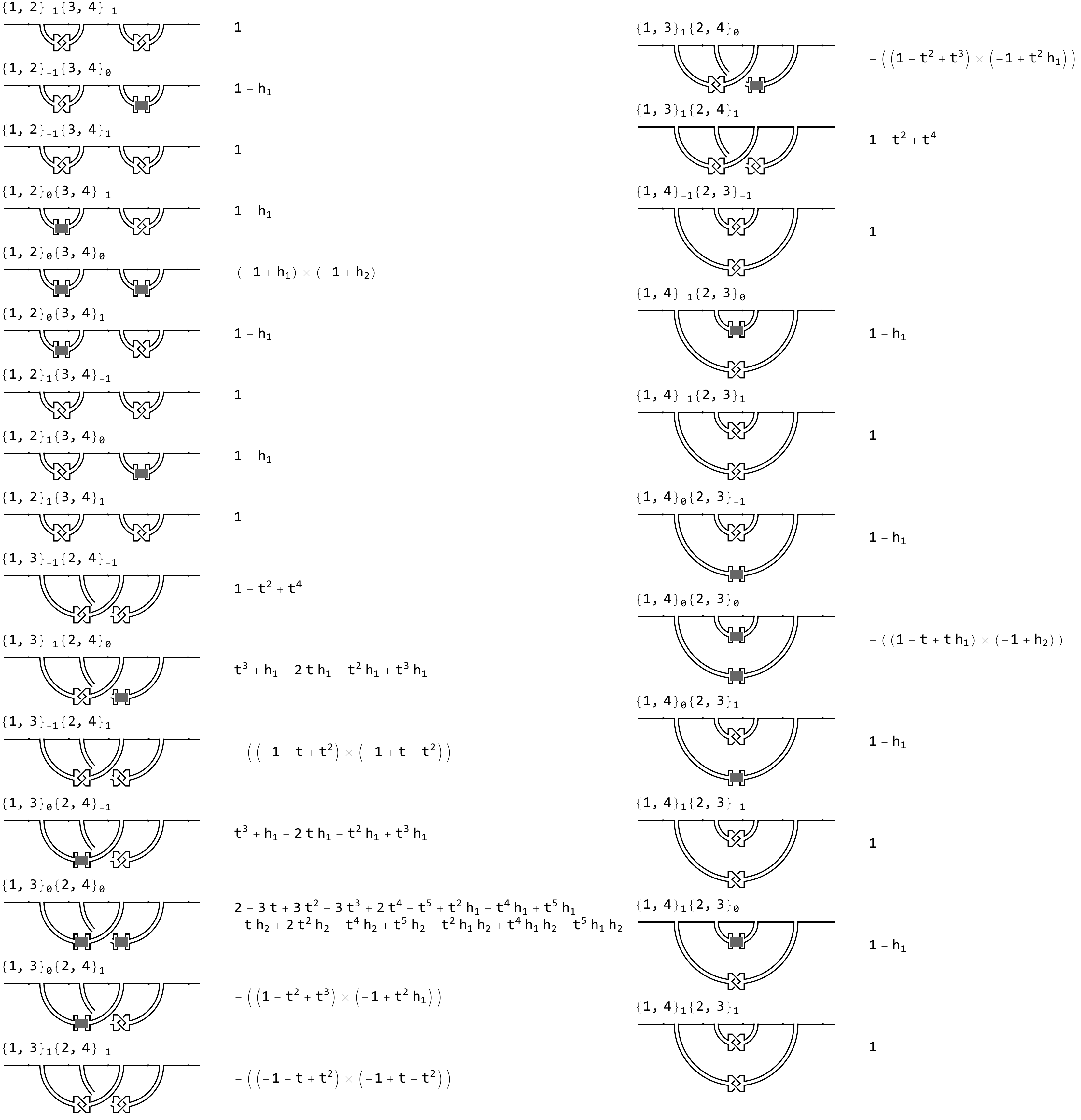}
\caption{All 27 generalized circuit topology diagrams with two contacts, with~their notation and value of $\DG_1$.}
\label{fig.TablegCT2}
\end{figure} 

\section{Conclusions and Further~Directions}

Generalized circuit topology studies folded linear molecular chains, taking into account both the way the chain bonds to itself (to form hard contacts) and the way it tangles in space (to form soft contacts or clasps). In~this study, we provided a theoretical basis for the language of circuit topology. For~such chains and their corresponding topological circuits, we furthermore introduced a new topological invariant $\DG$ that summarizes much of the essential topological features into
a single polynomial. A~simple computer implementation of the computation of the invariant was given in the Appendix~\ref{sec.Appendix} as well as a complete table
of the values of the invariant on all generalized circuit topology diagrams with at most three~contacts.

In future work it may be of interest to consider additional topological features such as rigid beads (vertices of valency two) and vertices of higher valency. The~alpha helices and beta sheets found in proteins would be an example of an application. Any such object may be modelled by some version of $\Gamma$-calculus. 
The general technique is always to first specify the appropriate Reidemeister moves, to~describe how our objects are supposed to move and interact. Next
we propose some general form with undetermined coefficients for our new objects. Applying $\Gamma$ to both sides of the Reidemeister moves 
yields algebraic equations for the coefficients that may be used to determine~them.

$\Gamma$ calculus is by no means the only knot invariant that could be extended to molecular chains. The~techniques used here are appropriate for a wide
class of knot invariants known as universal knot invariants~\cite{Ha06}. Roughly speaking there is an invariant for any ribbon algebra and the basic idea is to
place special algebra elements on the strands of the knot near each crossing and then multiply. 
The $\Gamma$-calculus is an example of such invariants for a particular algebra~\cite{BNV19,Vo18}
but depending on the task at hand other algebras may be more convenient. For~example if one would need a higher resolution than the one given by $\DG$ one good 
option is to choose the polynomial time knot polynomials of~\cite{BNV19}. So far these invariants have been introduced for tangles only but following the ideas of the present paper an extension to $H$-tangles seems possible. Unlike the more well known strong invariants, such as the Jones polynomial, Khovanov and Floer homology, the~polynomial time knot polynomials are still relatively efficient to compute,
 allowing for real world applications. In~future work we plan investigate such invariants for molecular~chains.
 
We have established a good way of distinguishing and manipulating folded linear molecular chains but many exciting challenges remain. Especially relevant for applications 
seems to be the design of a workable notion of relative topological distance between two given chains. This allows for mapping the topological evolution of biomolecular folds and the development of useful reaction coordinates for monitoring molecular folding~reactions.

%%%%%%%%%%%%%%%%%%%%%%%%%%%%%%%%%%%%%%%%%%
%\section{Patents}

% This section is not mandatory, but~may be added if there are patents resulting from the work reported in this~manuscript.

%%%%%%%%%%%%%%%%%%%%%%%%%%%%%%%%%%%%%%%%%%
\vspace{6pt} 

%%%%%%%%%%%%%%%%%%%%%%%%%%%%%%%%%%%%%%%%%%
%% optional
%\supplementary{The following are available online at \linksupplementary{s1}, Figure S1: title, Table S1: title, Video S1: title.}

% Only for the journal Methods and Protocols:
% If you wish to submit a video article, please do so with any other supplementary material.
% \supplementary{The following are available at \linksupplementary{s1}, Figure S1: title, Table S1: title, Video S1: title. A supporting video article is available at doi: link.} 

%%%%%%%%%%%%%%%%%%%%%%%%%%%%%%%%%%%%%%%%%%
\authorcontributions{Conceptualization, A.M. and R.v.d.V.; methodology, A.M. and R.v.d.V.; software, R.v.d.V.; investigation, A.M. and R.v.d.V.; resources, A.M.; writing---original draft preparation, A.M. and R.v.d.V.; writing---review and editing, A.M. and R.v.d.V. All authors have read and agreed to the published version of the~manuscript.}

\funding{\hl{~}%MDPI: Please add: ``This research received no external funding'' or ``This research was funded by NAME OF FUNDER grant number XXX.'' and  and ``The APC was funded by XXX''. Check carefully that the details given are accurate and use the standard spelling of funding agency names at \url{https://search.crossref.org/funding}, any errors may affect your future~funding.
}

\institutionalreview{\hl{~}%MDPI: In this section, please add the Institutional Review Board Statement and approval number for studies involving humans or animals. Please note that the Editorial Office might ask you for further information. Please add ``The study was conducted according to the guidelines of the Declaration of Helsinki, and approved by the Institutional Review Board (or Ethics Committee) of NAME OF INSTITUTE (protocol code XXX and date of approval).'' OR ``Ethical review and approval were waived for this study, due to REASON (please provide a detailed justification).'' OR ``Not applicable'' for studies not involving humans or animals. You might also choose to exclude this statement if the study did not involve humans or animals.
}

\informedconsent{\hl{~}%MDPI: Any research article describing a study involving humans should contain this statement. Please add ``Informed consent was obtained from all subjects involved in the study.'' OR ``Patient consent was waived due to REASON (please provide a detailed justification).'' OR ``Not applicable'' for studies not involving humans. You might also choose to exclude this statement if the study did not involve humans.

%Written informed consent for publication must be obtained from participating patients who can be identified (including by the patients themselves). Please state ``Written informed consent has been obtained from the patient(s) to publish this paper'' if applicable.
}

\dataavailability{Data is contained within the article or supplementary material.} 

%\acknowledgments{In this section you can acknowledge any support given which is not covered by the author contribution or funding sections. This may include administrative and technical support, or donations in kind (e.g., materials used for experiments).}

\conflictsofinterest{The authors declare no conflict of~interest.} 

%% Optional
%\sampleavailability{Samples of the compounds ... are available from the authors.}

%%%%%%%%%%%%%%%%%%%%%%%%%%%%%%%%%%%%%%%%%%
%% Only for journal Encyclopedia
%\entrylink{The Link to this entry published on the encyclopedia platform.}

%%%%%%%%%%%%%%%%%%%%%%%%%%%%%%%%%%%%%%%%%%
%% Optional
\abbreviations{Abbreviations}{
The following abbreviations are used in this manuscript:\\

\noindent 
\begin{tabular}{@{}ll}
CT & Circuit topology\\
gCT & Generalized circuit topology
\end{tabular}}

%%%%%%%%%%%%%%%%%%%%%%%%%%%%%%%%%%%%%%%%%%
%% Optional
\appendixtitles{yes} % Leave argument "no" if all appendix headings stay EMPTY (then no dot is printed after "Appendix A"). If~the appendix sections contain a heading then change the argument to "yes".
\appendixstart
\appendix
\section{Implementation and~Table}
\label{sec.Appendix}

In this appendix we present implementations of $\Gamma$ and $\DG$ in Mathematica~\cite{Wolfram}. We also give a few examples of how the code is used in practice. Finally a further table of gCT diagrams is provided. All the code can be downloaded from the second author's website
: \url{www.rolandvdv.nl/gCT}
(accessed September 2021)

\subsection{Implementation for $\Gamma$}

$\Gamma$ calculus is simple enough for the complete implementation to fit five lines of Mathematica code as shown below. 
The code remains quite close to the mathematics. A~minor deviation is that pairs $(\omega,A)$ are represented as $\Gamma[\omega,A]$.
The first line of code implements disjoint union. The~second line makes sure $\omega$ is always in its simplest form and $A$ should be written as a quadratic in $r,c$ with coefficients that are functions of $t$. The~remaining lines determine the merging and the value on the crossings. For~simplicity we write $X_{ij}$ instead of $\Gamma(X_{ij})$ and $X^{-1}$ is written more compactly as $\bar{X}$.
The effect of merging strands $i,j$ is written as $\mu_{i,j\to k}$.
\begin{figure}[H]
\includegraphics[width=0.99\linewidth]{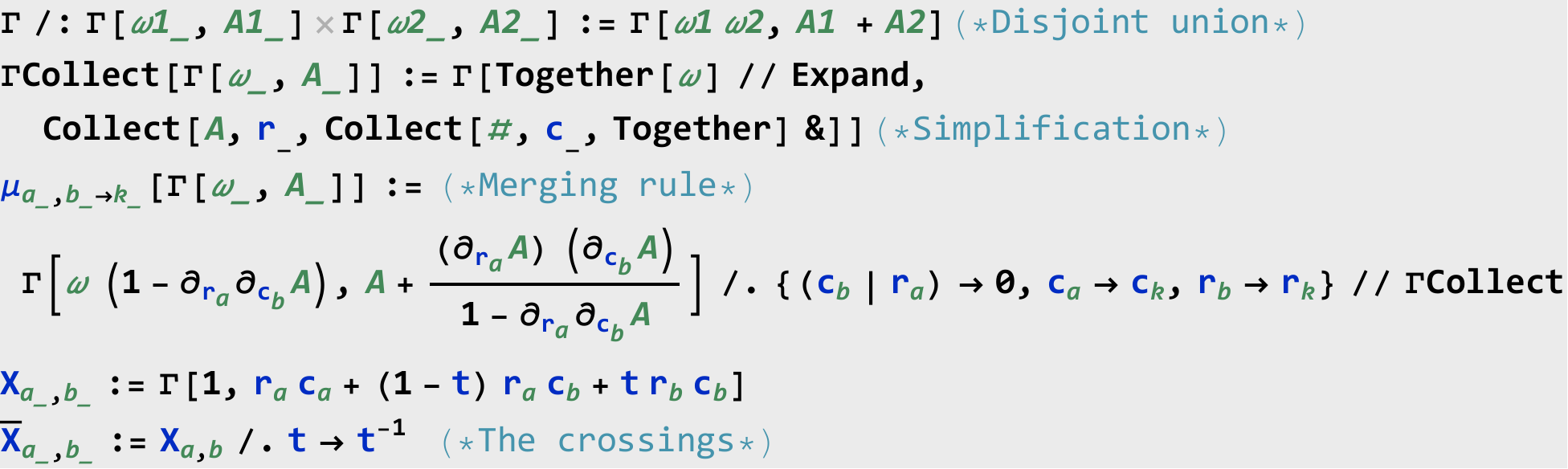}
\caption{A complete implementation of $\Gamma$-calculus.}
\label{fig.CodeGamma1}
\end{figure}

\vspace{-6pt}

\begin{figure}[H]
\includegraphics[width=10cm]{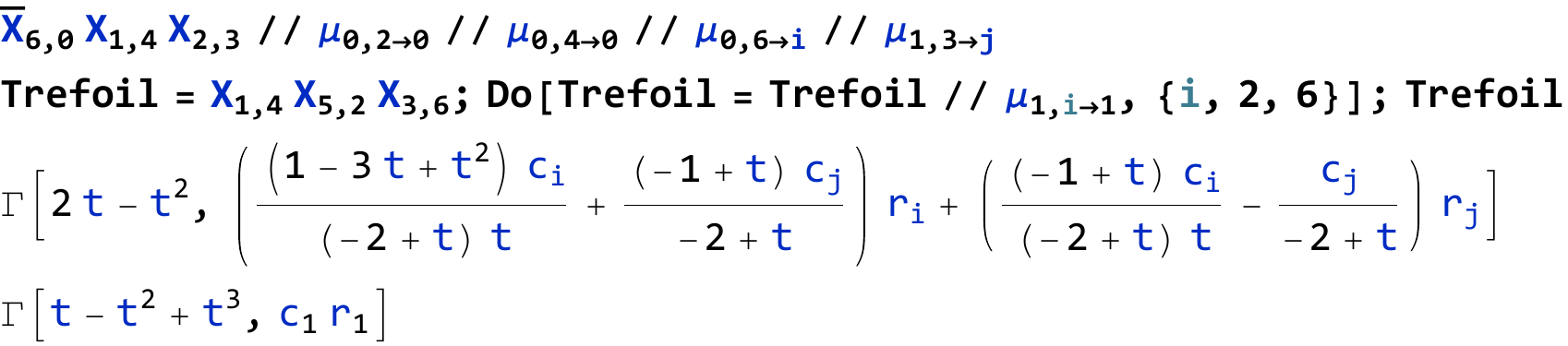}
\caption{Two examples of how to run the $\Gamma$ program.}
\label{fig.CodeGamma2}
\end{figure}

\subsection{Mathematica Implementation for $\DG$}

The implementation of $\DG$ is very similar to that of $\Gamma$ and builds on top of it. The~key difference is that
this time we require all strand labels to be positive in the sense that for each label $i$ the
label $-i$ is supposed to be distinct and not yet~occupied.

Pairs $(\omega,A)$ are still represented as $\Gamma[\omega,A]$. For~$\DG(X_{i,j})$ we write $XX_{ij}$ and $\overline{XX}_{ij}$ for the value of the negative crossing.
The value of the $H$-vertex $H^h_{ij}$ with label $h$ is denoted $HH_{ij}[h]$ and the merging is denoted $\mu\mu_{i,j\to k}$ 
Finally we add a small routine for comparing pairs $(\omega,A)$.
\begin{figure}[H]
\includegraphics[width=0.99\linewidth]{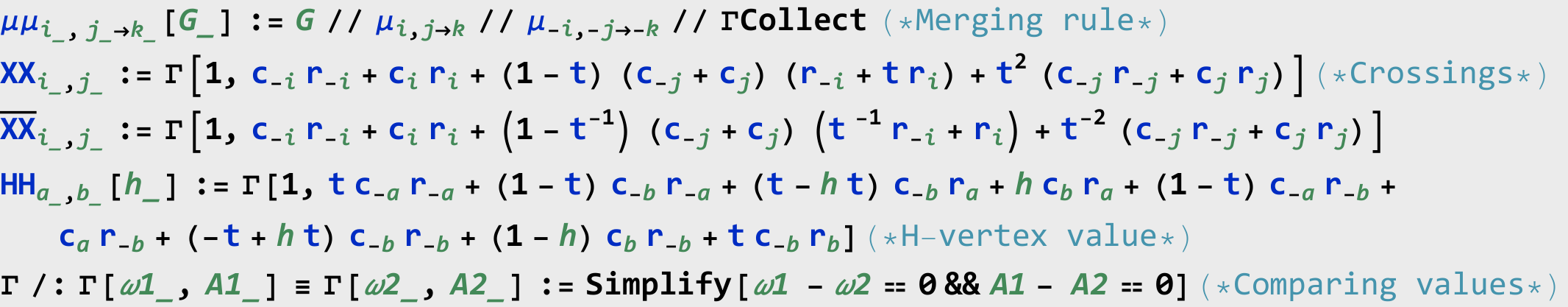}
\caption{A complete implementation of $\DG$-calculus.}
\label{fig.CodeDG1}
\end{figure}

To illustrate that the proof of invariance of $\DG$ is indeed automatic we provide the code to check
an instance of Reidemeister 3 and 4. The~other equations can (and should) be checked similarly by~computer.
\begin{figure}[H]
\includegraphics[width=12cm]{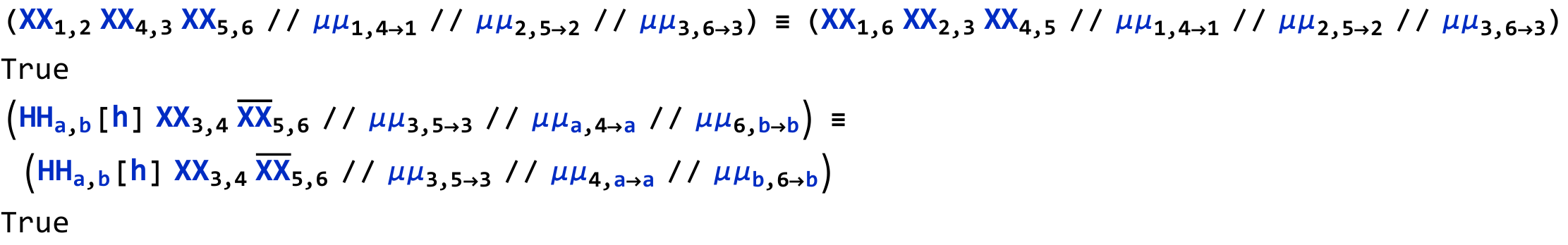}
\caption{Checking instances of Reidemeister 3 and~4.}
\label{fig.CodeDG2}
\end{figure}

\subsection{Tabulation of All gCT Diagrams with Three~Contacts}

Finally we list all gCT diagrams with three contacts, grouped by their value of $\DG_1$.
Trivial cases that contain a loose clasp $\{i,i+1\}_{\pm 1}$ have been ommitted as they are clearly equivalent to a gCT diagram with two contacts.
Although not perfect we see that $\DG$ is able to distinguish many types of chains. For~example, the~gCT diagrams with three hard contacts
and no soft contacts all get a different value of $\DG$.~With~a few exceptions such as the mirror image trefoils already seen in the two-contact table
most of the diagrams that have the same value of $\DG$ in this table are in fact equivalent.

\noindent
\includegraphics[width=.99\linewidth]{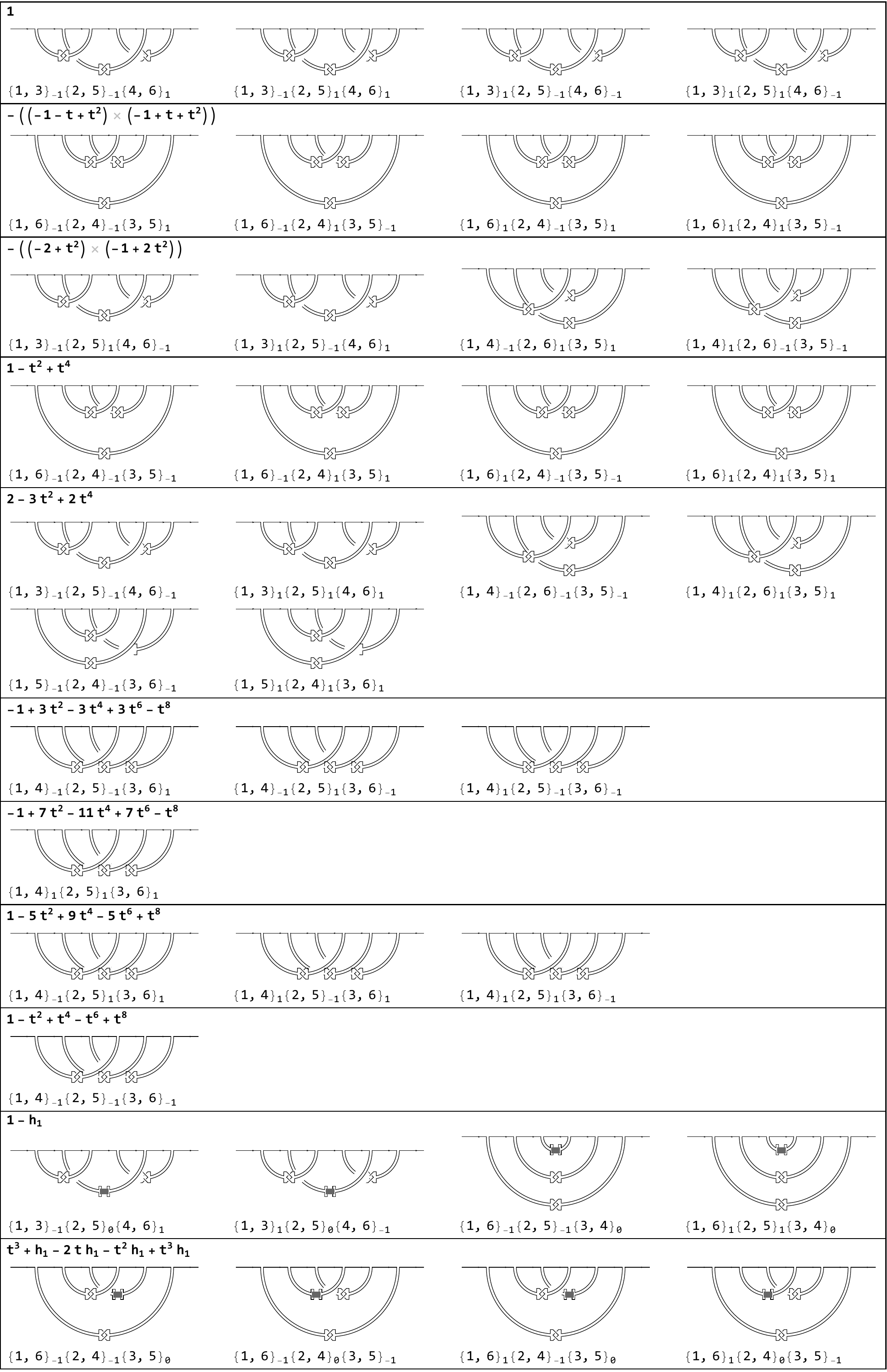}\\

\noindent
\includegraphics[width=.92\linewidth]{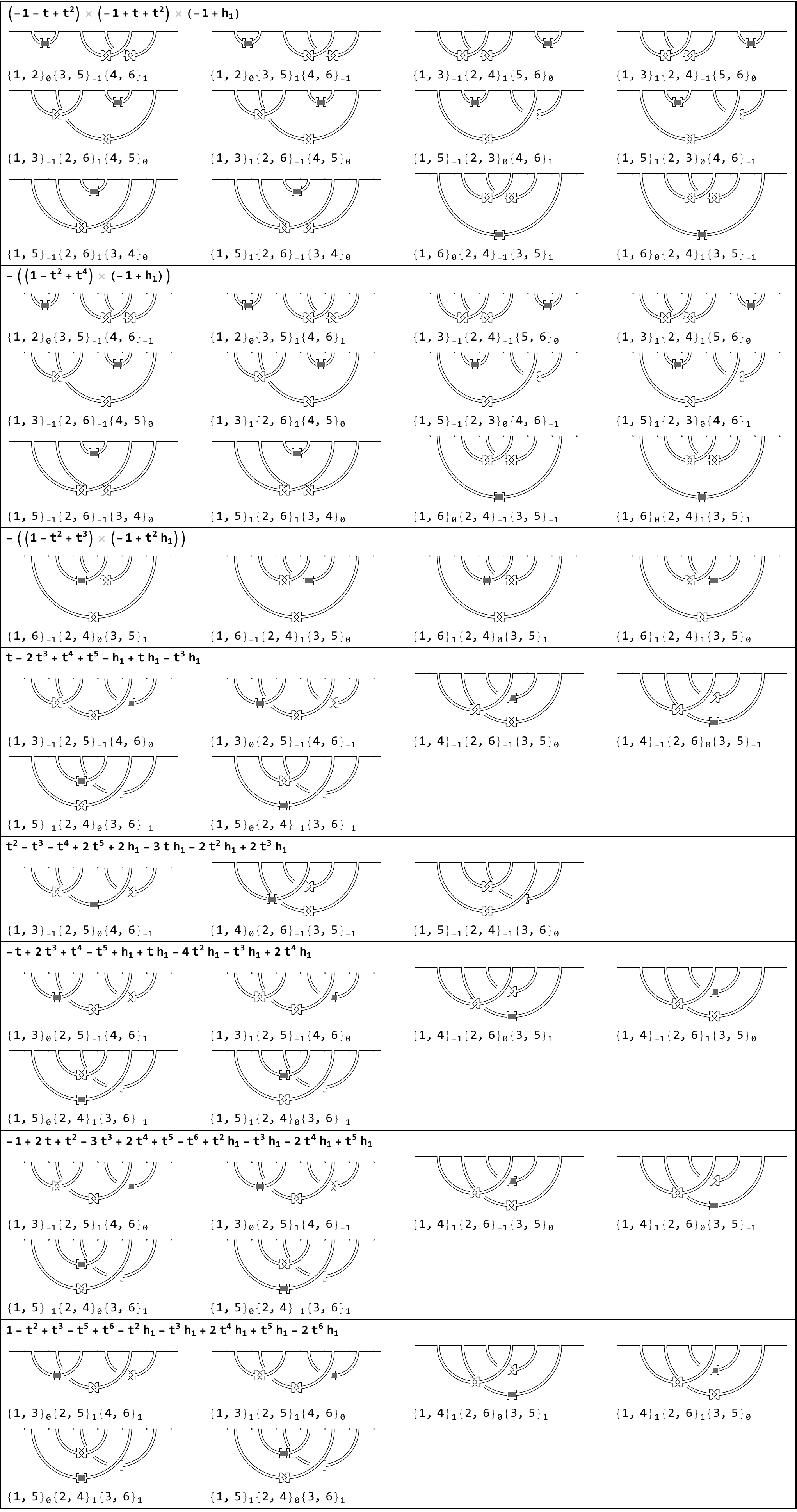}\\

\noindent
\includegraphics[width=.99\linewidth]{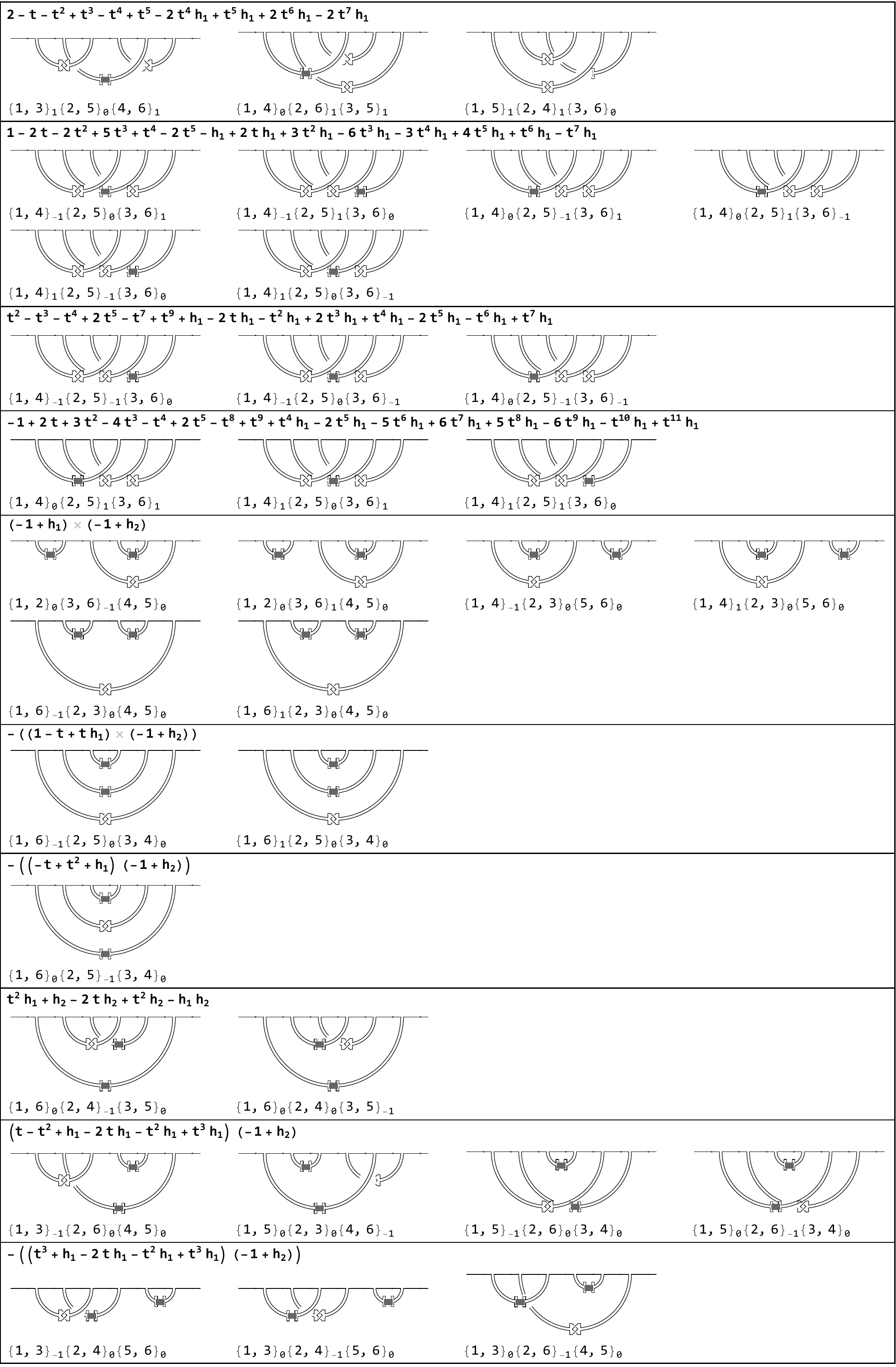}\\

\noindent
\includegraphics[width=.99\linewidth]{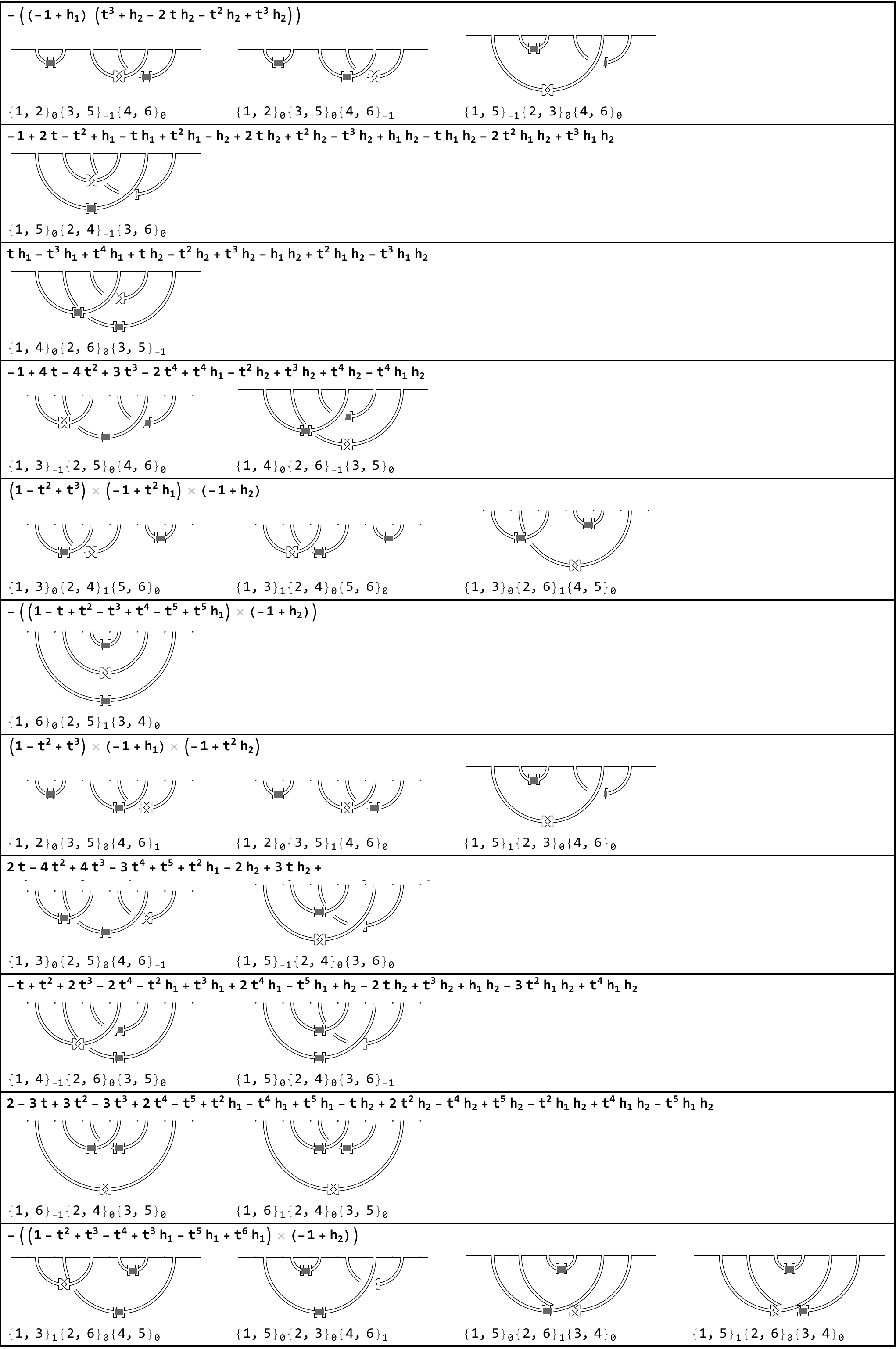}\\

\noindent
\includegraphics[width=.95\linewidth]{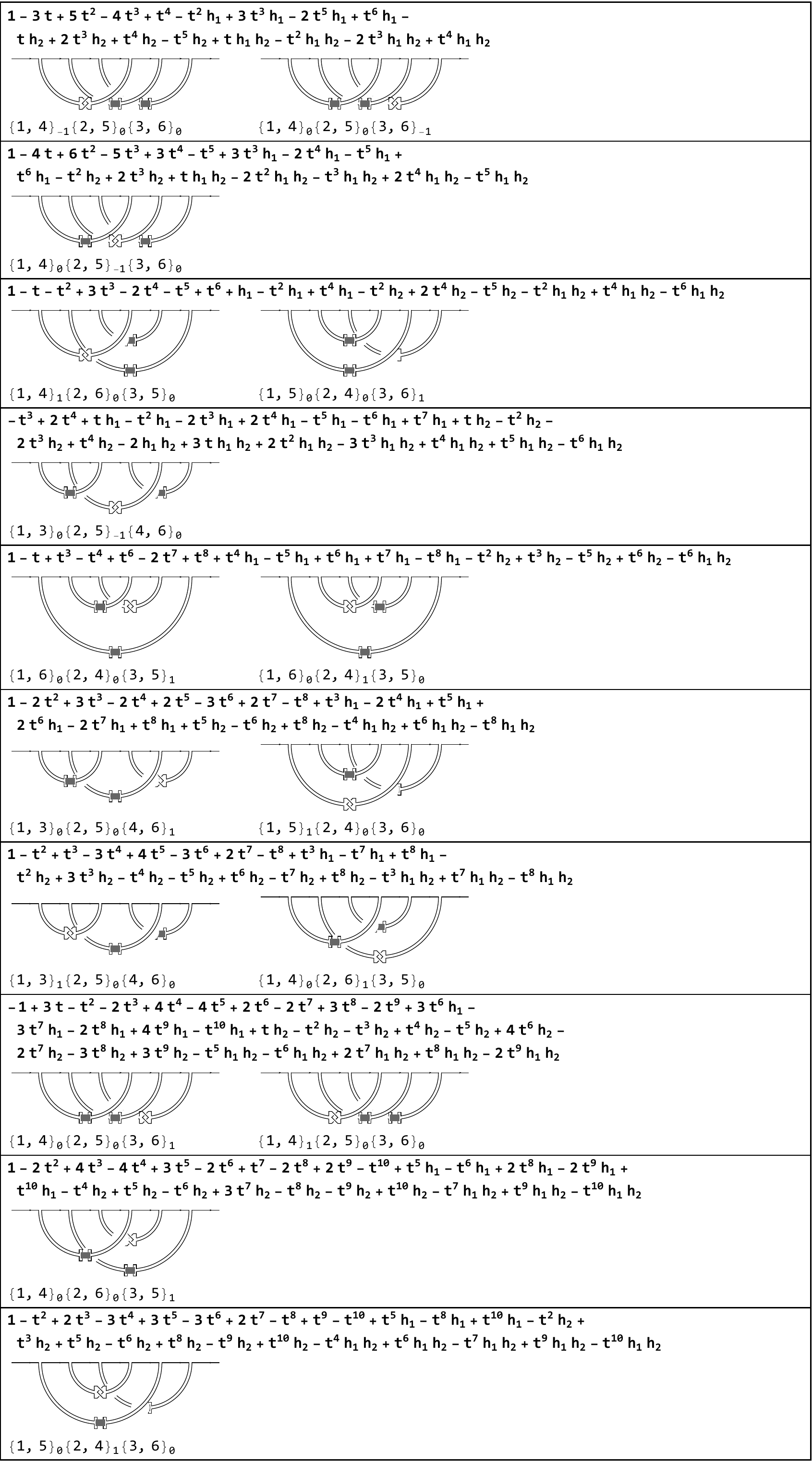}\\

\noindent
\includegraphics[width=.95\linewidth]{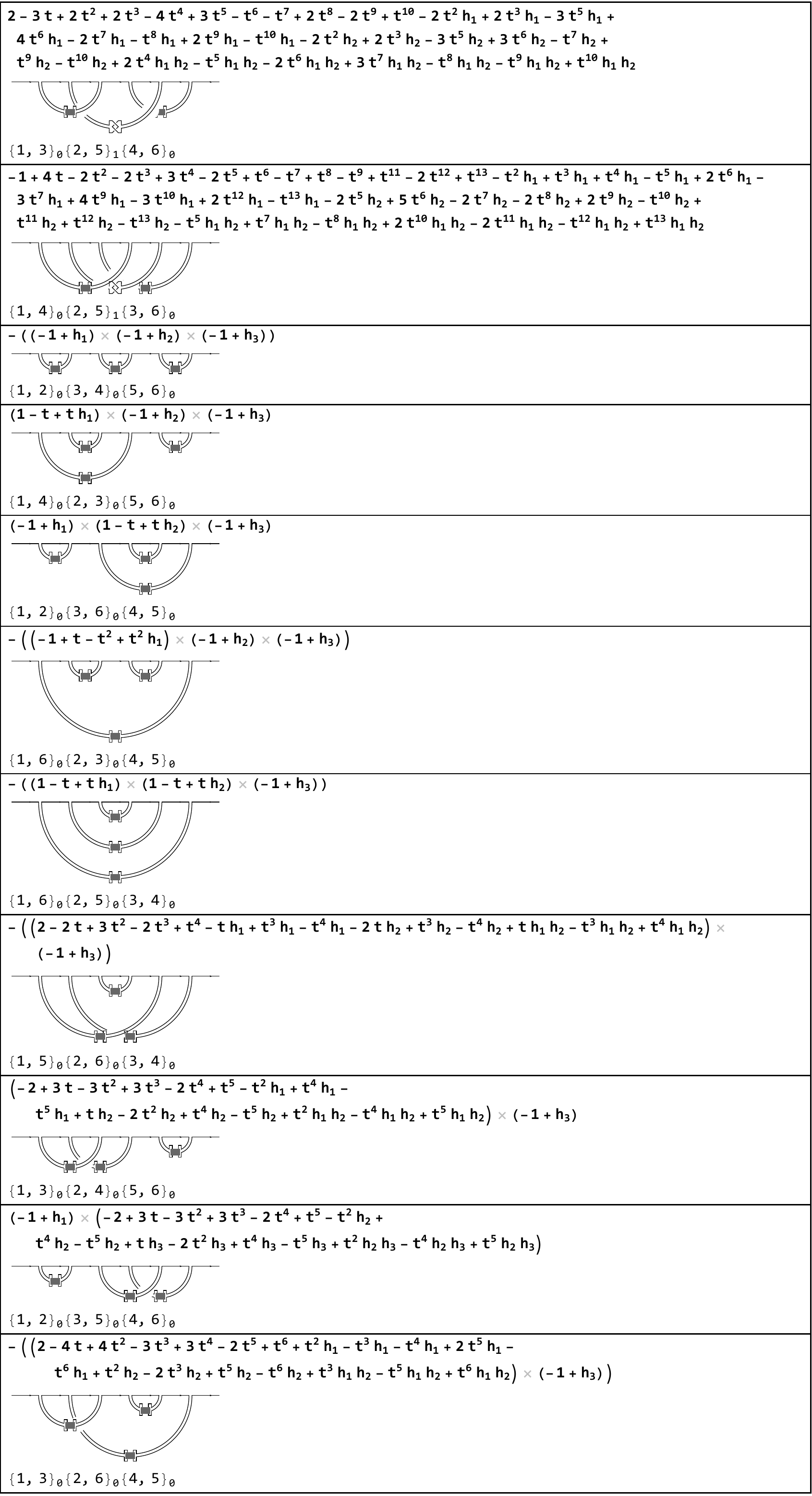}\\

\noindent
\includegraphics[width=.99\linewidth]{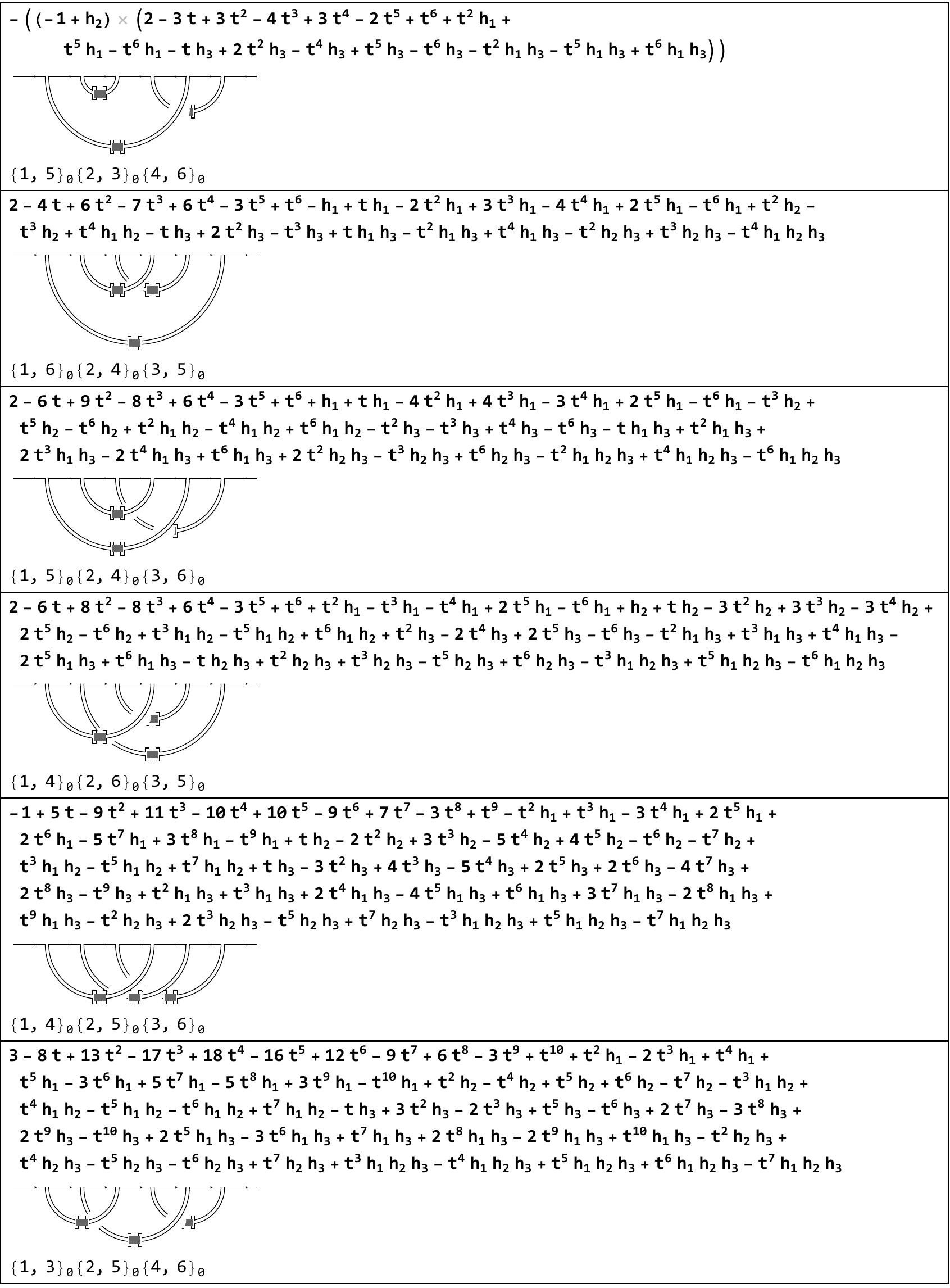}\\

\vspace{-12pt}

\end{paracol}

\reftitle{References}

\end{document}